\pdfoutput=1
\documentclass[12pt]{amsart}

\usepackage{amssymb,amsmath}
\usepackage{graphicx}
\usepackage{microtype}
\usepackage{mathdots}
\usepackage{amsbsy}
\usepackage{amscd}
\usepackage{amsthm}
\usepackage{mathrsfs}
\usepackage{verbatim}
\usepackage[table]{xcolor}
\usepackage[colorlinks,]{hyperref}
\usepackage{authblk}
\usepackage{mathtools}
\usepackage[english]{babel}
\usepackage[T1]{fontenc}
\usepackage[utf8]{inputenc}
\usepackage{tikz}
\usepackage{tikz-cd}
\usepackage{fullpage}
\usepackage{enumerate}
\usepackage[foot]{amsaddr}
\usepackage{comment}
\usepackage{tabularx}
\usepackage{tabulary}
\usepackage{setspace}

\definecolor{gry}{gray}{0.90}

\hypersetup{
 linkcolor = purple,
 citecolor = purple,
 filecolor = purple,
 urlcolor = black
 }

\theoremstyle{plain}
\newtheorem{Theorem}{Theorem}
\newtheorem{Corollary}[Theorem]{Corollary}
\newtheorem{Proposition}[Theorem]{Proposition}
\newtheorem{Lemma}[Theorem]{Lemma}

\theoremstyle{definition}
\newtheorem{Example}[Theorem]{Example}

\newtheorem{Definition}[Theorem]{Definition}
\newtheorem{DefProp}[Theorem]{Definition/Proposition}
\theoremstyle{remark}
\newtheorem{Remark}[Theorem]{Remark}

\newcommand{\Rnumeral}[1]{\uppercase \expandafter {\romannumeral #1 \relax}}

\numberwithin{Theorem}{section}

\newcommand{\N}{\mathbb{N}}
\newcommand{\RR}{\mathbb{R}}

\newcommand{\Z}{\mathbb{Z}}
\newcommand{\Q}{\mathbb{Q}}

\newcommand{\PP}{\mathbb{P}}
\newcommand{\G}{\mathbb{G}}

\newcommand{\mc}[1]{\mathcal{#1}} 
\newcommand{\ms}[1]{\mathscr{#1}} 
\newcommand{\mf}[1]{\mathfrak{#1}} 
\newcommand{\mb}[1]{\mathbb{#1}} 
\newcommand{\mt}[1]{\text{#1}}

\newcommand*{\defeq}{\mathrel{\vcenter{\baselineskip0.5ex \lineskiplimit0pt
                     \hbox{\scriptsize.}\hbox{\scriptsize.}}}=}

\newcommand{\Rnumeralum}[1]{\uppercase\expandafter{\romannumeral #1 \relax}}

\DeclareMathOperator{\Hom}{Hom}
\DeclareMathOperator{\Sh}{Sh}
\DeclareMathOperator{\Ker}{Ker}
\DeclareMathOperator{\res}{res}
\DeclareMathOperator{\rank}{rank}
\DeclareMathOperator{\Sym}{Sym}

\begin{document}

\title[Equivariant $K$-theory]{Equivariant $K$-theory of smooth projective spherical varieties}

\author[S.Banerjee]{Soumya Banerjee}
\address{Ben-Gurion University of the Negev, Israel}
\email{soumya@cs.bgu.ac.il}
\author[M.Can]{Mahir Bilen Can}
\address{Tulane University, New Orleans, USA}
\email{mcan@tulane.edu}

\keywords{Equivariant $K$-theory, equivariant Riemann-Roch, spherical varieties, minimal rank symmetric varieties, wonderful compactification}
\subjclass[2010]{19E08, 14M27}

\date{February 11, 2017}

\maketitle

\begin{abstract}
We present a description of the equivariant $K$-theory of a smooth projective spherical variety. This provides an integral
$K$-theory version of Brion's calculation of equivariant Chow-cohomology of such varieties. We consider the equivariant
$K$-theory of wonderful compactifications of minimal rank symmetric varieties. We obtain a formula for their structure
constants in terms of certain lower dimensional Schubert classes. This generalizes results of Uma on equivariant
compactifications of adjoint groups. 
\end{abstract}

%
%
%
%

\section{Introduction}
The foundations of equivariant algebraic geometry have matured enough to support a wealth of results that
match their topological counterparts. These results pave the way for studying equivariant (generalized) cohomology
theories of algebraic varieties with linear algebraic group actions. A broad class of such varieties: toric
varieties, and  spherical varieties also happen to admit concrete combinatorial descriptions. This leads to a rich
interaction between geometry,  and combinatorics.  In this paper we study the equivariant $K$-theory of smooth projective spherical varieties with a
particular emphasis on the wonderful compactifications of minimal-rank symmetric varieties.

\subsubsection{}
Building on the earlier ideas of Chang and Skjelbred, in their influential paper \cite{GKM98}, Goresky, Kottwitz and MacPherson introduced 
an effective method to calculate topological equivariant cohomology of an equivariantly formal space. 
The subsequent applications and refinements of their methods are now commonly referred to as ``GKM-type'' results. In
algebraic-geometry Brion obtained the first GKM-type presentation of equivariant rational Chow-ring of a smooth
spherical variety, see \cite{Brion97}, using the equivariant intersection theory of Edidin and Graham.

The study of equivariant algebraic $K$-theory for coherent sheaves was initiated in early 1980's by Thomason~\cite{Thomason83}\footnote{Note that in literature 
one denotes the Quillen $K$-functors on coherent sheaves by $G$ and on vector bundles by $K$. In this paper, we work with coherent sheaves 
and use $K$ unambiguously.}. Almost three decades after its inception, Vezzosi and Vistoli, in
\cite{vvi,VezzosiVistoli03}, substantially extended Thomason's work, which paves the way for GKM-type results in algebraic $K$-theory.

As the first application of Vezzosi and Vistoli results, Uma calculated presentations of the equivariant Grothendieck $K$-groups of the wonderful 
compactification of a semisimple reductive group of adjoint type, see \cite{Uma}. 

Working on a slightly different question, Joshua and Krishna, in \cite{JoshuaKrishna}, show that for a smooth projective spherical $G$-variety $X$ (defined over a field $k$), with a maximal torus $T \subset G$, 
there is a ring isomorphism 
\[ K_{T,0}(X) \otimes_{\Z} K_{T,\ast}(k) \cong K_{T,\ast}(X) .\] 
This result together with Uma's work recovers the complete equivariant $K$-theory in the group case. One of the goals of
our paper is to extend these ideas to other smooth and projective spherical varieties.

Recently, Anderson and Payne, in \cite{anp}, initiated a study of a operational bivariant theory associated to the Grothendieck groups of coherent
equivariant sheaves. Using this operational $K$-theory, Gonzalez extended Uma's results to possibly non-smooth spherical $G$-varieties admitting 
finitely many torus invariant curves (also called $T$-skeletal varieties); see \cite{gzz}. These results provide
answers for certain singular varieties, however at present they are not applicable to non-$T$-skeletal cases. Our work,
on the other hand, applies to all smooth projective spherical varieties. 

\subsubsection{} This paper was motivated by a question of Dan Edidin who drew our attention to extending the equivariant $K$-theory
computations in the group case to other spherical varieties. Our first result is an integral
$K$-theoretic analogue of Brion's presentation of equivariant Chow-groups for spherical varieties in \cite[Theorem 7.3]{Brion97}.

\begin{Theorem}\label{T:mt1}
Let $X$ be a smooth projective spherical $G$ variety. Let $T$ be a maximal torus of $G$, $\Phi_{G}$ the roots of $G$
with respect to $T$, $R(T)$ the representation ring of $T$ and $W$ the Weyl group of $G$. 
The $T$-equivariant $K$-theory $K_{T,*}(X)$ is the ring of ordered 
tuples $(f_{x})\in \prod_{x \in X^{T}} K_{*}(k) \otimes R(T)$ satisfying the following congruence conditions:
\begin{enumerate}

\item $f_{x} - f_{y} = 0 \mod (1-\chi)$ when $x,y$ are connected by a $T$-stable curve with weight $\chi$. 

\item $f_{x} - f_{y} = f_{x} -f_{z}  = 0 \mod (1-\chi)$ and $f_{y} - f_{z} = 0 \mod (1-\chi^{2})$, where $\chi \in
  \Phi_{G}$, and an irreducible component of the subvariety $X^{\Ker(\chi)}$ which contains the points $x,y,z$ is
  isomorphic to $\mb{P}^{2}$. There is an element in $W$ that fixes $x$ and permutes the point $y$ and $z$.

\item $f_{x} - f_{y} = f_{y} -f_{z}  = f_{z} - f_{w} = f_{x} - f_{w} = 0 \mod (1-\chi)$, where $\chi \in \Phi_{G}$, and an irreducible component of the subvariety $X^{\Ker(\chi)}$  which contains the points $x,y,z,w$  is isomorphic to $\mb{P}^{1}\times \mb{P}^{1}$. There is an element in $W$ that fixes two points and permutes the other two. 

\item $f_{x} - f_{y} = f_{z} -f_{w} = 0 \mod (1-\chi)$ and $ f_{y} - f_{z} = 0 \mod (1-\chi^{2n})$ and $f_{x} - f_{w} =0
  \mod (1-\chi^{n})$, where $n > 1$, and $\chi \in \Phi_{G}$. An irreducible component of the subvariety $X^{\Ker(\chi)}$  which contains the points $x,y,z,w$  is isomorphic to the ruled surface $\mb{F}_{n}$. There is an element in $W$ that fixes the points $x$ and $w$ and permutes $z$ and $y$. 

\end{enumerate}

The Weyl group $W$ acts on the torus fixed point set $X^{T}$ by permutation and the $G$-equivariant $K$-theory is obtained 
by taking $W$-invariants. 
\end{Theorem}

The natural strategy to prove this result is to reduce the problem to the computation of some small
dimensional (specifically rank one spherical $SL_{2}$ compactifications) spherical varieties. The classification of such
varieties is well known due to Ahiezer (see \cite{Ahiezer, Brion86}). The original problem is now tantamount
to calculating the equivariant $K$-theory of finitely many small dimensional varieties. At this point one encounters the
issue that not all spherical varieties are $T$-skeletal and can have positive dimensional families of torus invariant
curves which must be accounted for. A good example is the classical variety of complete quadrics. We circumvent this issue by rigidifying the 
problem using toric geometry. We recover the necessary $K$-theory from the toric case by using a co-base change 
theorem, see Proposition \ref{P:desc}.

When a spherical variety $X$  has finitely many  $T$-invariant curves, GKM-type theorems are much more tractable and
often contain additional combinatorial structures. The $T$-equivariant $K$-theory of flag varieties is perhaps the most
well studied example, see \cite{kok}. In this vein, a natural family of varieties are the wonderful
compactifications, due to De Concini and Procesi (see \cite{DP83}), of the minimal-rank\footnote{The
  minimal rank condition ensures that such varieties have finitely many torus fixed curves.} symmetric spaces. The work
of Tchoudjem, in \cite{Tchoudjem07}, gives a combinatorial description of the torus stable points and curves in the wonderful compactifications
of minimal rank symmetric varieties. When $G$ is a semisimple group of adjoint type a complete
classification of the irreducible components up-to isomorphism is known. These are the wonderful compactifications of
three families ({\romannumeral 1 }) $PSL_{2n}/PSp_{n}$, ({\romannumeral 2}) $PSO_{2n}/PSO_{2n-1}$, ({\romannumeral 3})
The group case : $G\simeq G\times G/ \Delta G$ (where $\Delta G$ is the diagonal embedding) and an isolated exceptional
case $E_{6}/F_{4}$. 
Many important results, in the study of equivariant cohomology, were obtained in the group case 
by Bifet, De Concini and Procesi in their seminal work in  \cite{BDP90}, and more recently, their has been extended by Strickland in~\cite{Strickland12}.

In the $K$-theory setting, the group case was thoroughly investigated bu Uma in \cite{Uma}. Building up on Tchoudjem's work we are able to
 generalize Uma's result to wonderful compactifications of all families of minimal-rank symmetric spaces.

\begin{Theorem} \label{T:2ndthm}[See Proposition~\ref{P:structure} for details]
Let $X$ denote the wonderful compactification of an irreducible minimal rank symmetric variety $G/H$. Then $K_{G,*}(X)$ has a decomposition of the form 
\[
K_{G,*}(X) = K_{S,*}(Y_0)\otimes R(T/S)^{W_H}, 
\]
where $Y_{0}$ is an affine toric variety, $W_H$ is the Weyl group of $H$ and $S$ is a maximal anisotropic subtorus of $T$.
\end{Theorem}

An immediate consequence of this theorem is that the $G$-equivariant $K$-theories of wonderful compactifications of $PSL_{n} \times PSL_{n}/ \Delta PSL_{n}$ and
$PSL_{2n}/PSp_{n}$ are identical. 

There is an important basis of the torus equivariant $K$-theory of flag variety, the Schubert basis, which has deep
combinatorial structure. In the minimal rank case we show that there is a basis which enjoys the same combinatorial
properties as the Schubert basis and the sturcture constants, in this basis, are related to that of an appropriate Schubert basis.

\subsubsection*{Outline of the paper}
We will present a brief outline of the contents of this paper.

In Section \ref{S:Kth}, we collect several results in equivariant $K$-theory that are crucial to the rest of the
paper. Some of these results are well known and some are new.

In Section \ref{S:Preliminaries}, we study the geometric structure of the fixed point locus $X^{S}$, where $S\subset T$ is a 
codimension one subtorus. This is an important step in reduction step that goes into the proof of Theorem \ref{T:mt1}.

In Section \ref{S:Kthsph}, we combine the results of the previous sections to prove Theorem \ref{T:mt1}. We start with
explicit presentation in the rank one case of the base cases and then bootstrap it to the general case. 

In Section \ref{S:applications}, we recall some relevant structure theory of wonderful compactifications of symmetric varieties 
of minimal rank. We present a proof of Theorem \ref{T:2ndthm} and several variants that extend the results of Uma, in
\cite{Uma}.
In the short appendix, we use equivariant Riemann-Roch theorem to relate our work with Brion's Chow-ring computations.

\subsection*{Acknowledgments} The present paper owes its existence to Dan Edidin who drew our attention to this question 
and generously answered our technical questions. We are extremely grateful to Michel Brion for his invaluable guidance, 
comments, corrections and encouragement in the final stages of this work. 

We thank Dave Anderson, Mikhail Kapranov, Sam Payne and Lex Renner for their crucial remarks and suggestions during various 
stages of this work. The first author would like to thank Amnon Besser, Ilya Tyomkin and Amnon Yekutieli for suggestions and answers to technical questions.

\subsubsection{Assumptions/ Notations}\label{sbsb:assump}
A variety is a reduced scheme of finite type defined over an algebraically closed field $k$ of characteristic zero. It
is allowed to have irreducible components. Points of a variety will always mean closed points. The representable functor, from schemes over $k$ to the category of sets, associated to a scheme $X$ will be denoted by $\underline{X}$. The set $\underline{X}(A)$ will denote the $A$ valued points of $X$ for any $k$-algebra $A$. 
  
Unless otherwise stated, the linear group $G$ is always connected and reductive.  
A $G$-variety $X$ is a normal variety with an algebraic action of $G$. Given any closed subgroup $H \subset G$ the
neutral component of $H$ will mean the component of $H$ containing the identity. The semisimple part of a reductive
group $G$ will be denoted by $G^{ss}$.  

We recall that \emph{rank} has two, potentially confusing, meanings. The rank of a linear groups is its semisimple-rank
where as the rank of a spherical variety of homogeneous space is defined in Section \ref{S:fixed point loci}. 

All tensor products will be over $\Z$ unless specified.

%
%
%
%

\section{Results in $K$-theory}\label{S:Kth}
In this paper the term \emph{equivariant algebraic $K$-theory} will mean the algebraic $K$-theory of the category of
equivariant coherent sheaves on a $G$-variety $X$. The foundational results were obtained by Thomason; see \cite{Thomason83}. Informally speaking, in this
theory, one applies Quillen's $Q$-construction to the abelian category of equivariant sheaves.

\begin{Definition}[Equivariant sheaf]
Let $a$ denote the action map $a: G\times X\rightarrow X$ and let $p_X: G\times X\rightarrow X$ denote the second projection. 
A $G$-equivariant sheaf on $X$ is a pair $(\mc{F}, \phi)$ where $\mathcal{F}$ is a coherent sheaf on $X$ and  $\phi$ is an isomorphism $\phi: a^{\ast}\mathcal{F} \rightarrow p_{X}^{\ast} \mathcal{F}$ of sheaves on $G \times X$, which satisfies a natural cocycle condition on $G\times G \times X \rightrightarrows G\times X $. 
\end{Definition}

A morphism of equivariant sheaves, $(\mc{F}, \phi)$ and $(\mc{F}^{\prime}, \phi^{\prime})$, is a morphism of sheaves which commutes with isomorphisms $\phi$ and $\phi^{\prime}$.
We denote the category of $G$-equivariant sheaves on $X$ by $\Sh_{G} (X)$. The equivariant $K$-groups, denoted by $K_{G,\ast}(X)$ are defined as the
homotopy groups of the loop-space of the nerve of $QSh_{G}(X)$ (where $QSh_{G}(X)$ is the Quillen $Q$ construction
applied to $\Sh_{G}(X)$). One can similarly define $K_{H,\ast}(X)$ for any closed subgroup $H \subset G$.  

These $K$-groups admit several homomorphisms resulting functorially from equivariant maps between spaces and homomorphisms of the
structure group $G$ acting compatibly on a fixed space. As we will see, this is often very useful for calculations.
We start with a general result that relates, for a group $G$ and a closed subgroup $H$, the $G$-equivariant $K$-theory
to $H$-equivariant $K$-theory. 

\begin{Proposition}[Faddeev-Shapiro Lemma, \cite{Merkurjev97}]\label{Faddeev-Shapiro}
Let $H$ be a closed subgroup of the algebraic group $G$. 
Then for any $G$-variety $X$, the inclusion map $X \hookrightarrow X \times G/H$ 
defined by $x\mapsto (x, eH)$ induces the isomorphisms 
\begin{eqnarray*}
K_{G,n} (X\times G/H) \simeq K_{H,n}(X), 
\end{eqnarray*}
for all $n\geq 0$.
\end{Proposition}

Notice that we do not assume $G$ is reductive in the above proposition. When $G$ is a reductive group and $B \subset G$ is a Borel subgroup there is an useful refinement. 

\begin{Proposition} [Merkurjev, \cite{Merkurjev97}]
Let $G$ be a (split) reductive group and $B$ a Borel subgroup and $X$ be a smooth projective $G$-variety then the natural map 
\[
\theta : R(B) \otimes_{R(G)} K_{G,n}(X) \rightarrow K_{B,n}(X)
\] 
is an isomorphism for all $n \geq 0$.
\end{Proposition}

Using the structure theory of solvable groups, we can decompose the Borel group $B$ into a maximal torus $T$
and a unipotent subgroup $U$. As an algebraic variety $U$ is isomorphic to an affine space. The homotopy-invariance
property of equivariant $K$-theory then implies 
$K_{B,n} (X) \simeq K_{B,n} (X\times B/T) \simeq K_{T,n}(X)$; see \cite[Theorem 4.1]{Thomason83} for details.

\begin{Remark}\label{R:Gequi}
In particular, we have 
\[K_{T}(X) = R(T) \otimes_{R(G)} K_{G}(X).\]
\label{R:SteinbergPittie}
When $G$ is simply connected, a result of Steinberg (see \cite{Steinberg75}) shows that
\begin{align*}
R(G) \simeq R(T)^W  \text{ and }  R(T) \simeq R(G)\otimes \Z[W],
\end{align*}
where $W$ is the Weyl group of $(G,T)$. The first isomorphism is an isomorphism of rings and the latter is only an isomorphism of $R(G)$-modules with a
compatible action of $W$. Consequently, we recover $K_{G}(X)$ as the space of $W$-invariants in $K_{T}(X)$.
\end{Remark}
In the next Proposition we consider a refinement of Proposition
\ref{Faddeev-Shapiro} for torus actions.

\begin{Proposition} \label{L:Edidin's tip}
Let $T$ be an algebraic torus and $T^{\prime} \subset T$ is a fixed codimension one subtorus. Let $X$ be a projective $T$-variety. Then the canonical map 
\[
R(T^{\prime}) \otimes_{R(T)} K_{T,n}(X) \rightarrow K_{T^{\prime},n}(X) 
\] is an isomorphism.
\end{Proposition}

\begin{proof}
Let $\mathcal{L}$ be a one-dimensional representation of $T$ such that $\mathcal{L}\setminus \{ 0\} = T/T^{\prime}$. 
Let $j$ denote the obvious map $j: X \times T/T^{\prime} \hookrightarrow X \times \mathcal{L}$, 
and $p: X\times \mathcal{L} \rightarrow X$ denote the projection map. 
The pullback map $p^{\ast}: K_{T,n}(X) \rightarrow K_{T,n}(X \times \mathcal{L}) $ is an isomorphism by the homotopy
invariance property. Let $\res:\Sh_{T}(X) \rightarrow \Sh_{T^{\prime}}(X)$ denote the canonical restriction map. In the
commutative triangle in eqn.(\ref{E:res}) the vertical equality is a consequence of Proposition \ref{Faddeev-Shapiro}.
 
\begin{equation} \label{E:res}
\begin{tikzcd}
K_{T,n}(X) \arrow{rd}{\text{res}} \arrow{r}{(pj)^{\ast}} & K_{T,n}(X \times T/T^{\prime}) \arrow[equal]{d} \\
& K_{T^{\prime},n}(X).
\end{tikzcd}
\end{equation}

Now consider the decomposition of the total space $X \times \mc{L}$ 
\[
\begin{tikzcd}
 X= X \times \{0\} 
 \arrow[hookrightarrow]{r}{i}& X \times \mathcal{L}  \arrow[hookleftarrow]{r}{j}&  X \times T/T^{\prime} 
 \end{tikzcd} 
\]
into a closed set $X \times \{0\}$ and its open complement. The terms of the localization long-exact sequence in
equivariant $K$-theory fit into the commutative diagram below. 
\begin{equation}\label{E:ressplit}
\begin{tikzcd}
{ }   & K_{T,n}(X) \arrow[dashed]{r}{i_{\ast}p^{\ast-1}}\arrow[equal]{d}  
& K_{T,n}(X) \arrow{r}{\res} \arrow{d}{p^{\ast}} & K_{T^{\prime},n}(X)  \arrow[equal]{d} & {}  \\
 \ldots \arrow{r} &  K_{T,n}(X) \arrow{r}{i_{\ast}} &  K_{T,n}(X\times\mathcal{L}) \arrow{r}{j^{\ast}}  
 & K_{T,n}(X \times T/T^{\prime})  \arrow{r}  & \ldots  
\end{tikzcd}
\end{equation}

Thanks to \cite [\Rnumeral{5}, Corollary 27]{Handbook}, the top horizontal row in Diagram (\ref{E:ressplit}) is split
exact and the kernel of $\res$ is generated by the ideal $(1-\chi)$ in $K_{T,n}(X)$; where $\chi$ is defined by the short exact sequence \[ 1\rightarrow T^{\prime} \rightarrow T \stackrel{\chi}{\rightarrow} \mathbb{G}_{m} \rightarrow 1.\]

As a result, we have \[ K_{T^{\prime},n}(X) = K_{T,n}(X)/ (1-\chi)K_{T,n}(X) = R(T^{\prime}) \otimes_{R(T)} K_{T,n}(X).\]
\end{proof}

\begin{Corollary}
Let $T$ be an algebraic torus and $T^{\prime}$ is any subtorus. Let $X$ be any projective $T$-variety then the canonical map $R(T^{\prime}) \otimes_{R(T)} K_{T,n}(X) \rightarrow K_{T^{\prime},n}(X)$  is an isomorphism. 
\end{Corollary}
\begin{proof}[Sketch of the proof]
We use induction on the codimension of $T^{\prime}$ in the torus $T$. We can and choose a chain of subtori 
\[
T'=T_0 \subset T_1 \subset \cdots \subset T_r=T
\] 
such that each $T_{i}$ is codimension one in $T_{i+1}$.
The assertion follows from using Proposition \ref{L:Edidin's tip} inductively at each step.
\end{proof}

\begin{Proposition}\label{P:desc}
Consider a fixed subtorus $T^{\prime} \subset T $. Let $X$ be a $T$-variety such that $T^{\prime}$ acts trivially on $X$. 
Then we have the following formula \[K_{T,*}(X) \cong R(T) \otimes_{R(T/T^{\prime})} K_{T/T^{\prime}, *}(X).\] 
\end{Proposition}

\begin{proof}
The torus $T/T^{\prime}$ acts on $X$ and the groups $K_{T/T^{\prime}, *}(X)$ are defined. Consider the following diagram
\[
\begin{tikzcd}
 {} & T/T^{\prime} \times X \arrow[shift left= 0.75ex]{d}[swap]{\overline{a} \;\;\;\smallskip} \arrow[shift right= 0.5 ex]{d}{   \;\;\; \overline{pr}_{X}  }  &   \\
 T \times X \arrow{ru}{\pi}  \arrow[shift left= 0.75ex]{r}{a} \arrow[shift right= 0.75ex]{r}[swap]{pr_{X}} & X 
  \end{tikzcd}
\]
where $a$ is the action map and $\overline{a}$ is the induced action map. Let $pr_{T}$ (resp. $pr_{X}$) denote the
projection map from $T \times X$ to $T$ (resp. $X$) and $\overline{pr}_{X}$ denotes the projection to $X$. We will
identify the isomorphisms
 $a^{\ast} \cong \pi^{\ast} \overline{a}^{\ast}$  and $pr_{X}^{\ast} \cong
\pi^{\ast} \overline{pr}_{X}^{\ast}$ on sheaves.

Let $\chi$ be any character of $T^{\prime}$. We fix an isomorphism $\phi: T^{\prime} \times T/T^{\prime} \rightarrow
T$ throughout
which induces the given embedding $T^{\prime} \subset T$. Let $\widehat{\chi}$ denote the unique character of $T$ which
restricts to $\chi$ on $T^{\prime}$ is trivial on $T/T^{\prime}$ (via $\phi$). We introduce the notion of
``$\chi$-twist'' of an equivariant sheaf $(\mc{F}, \mu) \in \Sh_{T/T^{\prime}}(X)$ (which depends on $\widehat{\chi}$). 
To this end, consider the map 
\begin{equation}\label{E:twst}
\chi \cdot \mu \defeq pr_{T}^{\ast}(\widehat{\chi}) \cdot \pi^{\ast}(\mu) : a^{\ast}(\mc{F}) \rightarrow pr_{X}^{\ast}(\mc{F})
\end{equation}
between sheaves $a^{\ast}(\mc{F})$ and $pr_{X}^{\ast}(\mc{F})$ on $T \times X$.
Given any character $\eta$ of $T$ the regular function $pr^{\ast}_{T}(\eta)$ is an invertible function on $T \times
X$. So the isomorphism $\chi \cdot \mu$ in the above equation \eqref{E:twst} is an isomorphism of sheaves and $(\mc{F}, \chi \cdot \mu)$ is a $T$-equivariant sheaf on $X$. 

The following is immediate from the definitions. 
\begin{align}\label{E:orth}
 \Hom_{\Sh_{T}(X)} ((\mc{F}, \chi \cdot \mu) , (\mc{G}, \chi^{\prime} \cdot \lambda)) = 
  \begin{cases} 
   0  & \text{if } \chi \neq \chi^{\prime} \\
   \Hom_{\Sh_{T/T^{\prime}}(X)}((\mc{F}, \mu), (\mc{G}, \lambda))       & \text{ if } \chi = \chi^{\prime}
   \end{cases}
\end{align}

We have thus constructed a functor $\mc{T}_{\chi}: \Sh_{T/T^{\prime}}(X) \rightarrow \Sh_{T}(X)$ which is full and
faithful. Given any object $\mc{F}$ of $\Sh_{T/T^{\prime}}(X)$ we call $\mc{T}_{\chi}(\mc{F})$ the $\chi$-twist of
$\mc{F}$ and let $\Sh_{T}(X)_{\chi}$ denote the subcategory $\mc{T}_{\chi}(\Sh_{T/T^{\prime}}(X))$ of
$\Sh_{T}(X)$ and for any character $\chi^{\prime}$ of $T^{\prime}$ (also extended to $T$ via $\phi$) the functor $\mc{T}_{\chi^{\prime}}$ defines equivalence of categories
\begin{equation}\label{E:sheq} 
\mc{T}_{\chi^{\prime}}: \Sh_{T}(X)_{\chi}  \stackrel{\cong}{\longrightarrow} \Sh_{T}(X)_{\chi^{\prime}\cdot\chi} = \mc{T}_{\chi^{\prime} \cdot \chi}(\Sh_{T/T^{\prime}}(X)).
\end{equation} 

The map of group schemes $T \rightarrow T/T^{\prime}$ is faithfully flat, by  \cite[Theorem 14.1]{Waterhouse}, so the
functor $\pi^{\ast}$ induces an equivalence between the category of sheaves $\Sh(X \times T/T^{\prime})$ and $\Sh(X
\times T)$ (see \cite[Chapter 6]{neron}). Let $\mc{D}$ denote the inverse equivalence. Then given any equivariant sheaf
$(\mc{F}, \alpha)$ in $\Sh_{T}(X)$ we get an isomorphism \[\mc{D}(\alpha): \overline{a}^{\ast}(\mc{F}) \rightarrow
  \overline{pr}_{X}^{\ast}(\mc{F})\] of sheaves on $X \times T/T^{\prime}$. This makes $(\mc{F}, \mc{D}(\alpha))$ into a
$T/T^{\prime}$-equivariant sheaf. The functor $\pi^{\ast}$ also preserves equivariant sheaves. This follows from the
observation that $\pi^{\ast}$ is identical to the functor $\mc{T}_{\chi_{e}}$, where $\chi_{e}$ is the trivial character.

Consider the abelian sub-category $\mc{C}$ of $\Sh_{T}(X)$ whose objects are $T$-equivariant sheaves
$(\mc{F}, \chi \cdot \mc{D}(\alpha))$, for any character $\chi$ of $T^{\prime}$, and morphisms \[ \Hom_{\mc{C}}((\mc{F},
  \chi \cdot \mc{D}(\alpha)), (\mc{G}, \chi^{\prime} \cdot \mc{D}(\alpha))) \defeq \Hom_{\Sh_{T}(X)}((\mc{F}, \chi \cdot
  \mc{D}(\alpha)), (\mc{G}, \chi^{\prime} \cdot \mc{D}(\alpha)))\] (see eqn. (\ref{E:orth}) above). Then $\mc{C}$ is equivalent to $\Sh_{T}(X)$ since any object $(\mc{F}, \mu)$ in $\Sh_{T}(X)$ is of the form $(\mc{F},
\chi_{e} \cdot \mc{D}(\mu))$ and the Hom-sets in $\Sh_{T}(X)$ are graded by the lattice of characters of $T^{\prime}$. The $K$-theory of $\Sh_{T}(X)$ is the same as the $K$-theory of $\mc{C}$.  

It follows from Quillen's Q-construction, that the space $BQ(\mc{C})$ is a disjoint discrete $\Hom(T^{\prime},
\mb{G}_{m})$-fold covering space of $BQ(\Sh_{T/T^{\prime}}(X)).$ The functors $\mc{T}_{\chi}$ induce an action of
$R(T^{\prime})$  on $BQ(\mc{C})$ (see eqn. \ref{E:sheq}). Calculating the homotopy groups we get \[ K_{T,*}(X) =
  K_{T/T^{\prime},*}(X) \otimes_{K_{*}(k)} R(T^{\prime}).\]

The module $K_{T/T^{\prime},*}(X)$ is an $R(T/T^{\prime})$ module and $R(T) = R(T^{\prime}) \otimes R(T/T^{\prime})$. This proves the assertion.  
  \end{proof}

\begin{Remark}
We note an alternate way to prove Proposition \ref{P:desc} without explicitly using the Q-construction: Theorem 1.2 of \cite{JoshuaKrishna} reduces the original
problem to understanding the relation between the Grothendieck $K$-groups $K_{T,0}(X)$ and $K_{T/T^{\prime},0}(X)$. This
follows from the constructions in the first part of the above proof. 
\end{Remark}

\subsection{Equivariant $K$-theory of toric varieties} \label{SS:toric}
We recall the fundamental result of Vezzosi and Vistoli which is at the foundation of most approaches to equivariant
algebraic $K$-theory of $G$-varieties. 

\begin{Theorem}[Vezzosi-Vistoli \cite{VezzosiVistoli03}]\label{VV}
Suppose $G$ is a diagonalizable group acting on a smooth proper scheme $X$ defined 
over a perfect field; denote by $T$ the toral component of $G$, that is the maximal 
subtorus contained in $G$. Then the restriction homomorphism on $K$-groups 
$K_{G,*}(X) \rightarrow K_{G,*}(X^{T})$ is injective, and its image equals the intersection 
of all images of the restriction homomorphisms $K_{G,*}(X^S) \rightarrow K_{G,*}(X^{T})$ 
for all subtori $S\subset T$ of codimension 1. 
\end{Theorem}

Among its many applications, the theorem provides a complete description of torus equivariant $K$-theory of smooth toric
varieties. Toric varieties will play an important role in  the following sections so we recall some basic results from the theory
of toric varieties. We refer the reader to \cite{Fulton} for details.

Let $T$ be an algebraic torus and we denote its lattice of characters (respectively, co-characters) by $M_{T}$
(respectively, by $N_{T}$). Recall that $R(T) = \Z[M_{T}]$ and $T = \Hom_{gp}(M_{T}, \mathbb{G}_{m})$. Given any
homomorphism   $\phi: T^{\prime} \rightarrow T$ we have a map $\phi^{\ast}: \Z[M_{T}] \rightarrow \Z[M_{T^{\prime}}]$
which makes $\Z[M_{T^{\prime}}]$ a $\Z[M_{T}]$ module. When $\phi$ is a closed embedding the induced map $\phi^{\ast}:
M_{T} \rightarrow M_{T^{\prime}}$ is surjective. We associate the subgroup $M_{T^{\prime}}^{\perp} \subset M_{T}$ which
consists of all characters of $T$ that are trivial on $T^{\prime}$ and we have a non-canonical splitting $M_{T}= M_{T^{\prime}} \oplus M_{T^{\prime}}^{\perp}$. When the underlying torus is clear from context we will drop the subscript $T$ from the (co-)character lattices.

A rational fan $\Delta$ in $N_{\RR} \defeq N \otimes \RR$ defines a toric variety,  denoted by $X(\Delta)$, with a dense
open set $T = \Hom_{gp}(M, \mathbb{G}_{m})$. We will exclusively work with rational fans and simply refer to
them as fans. Let $\Delta_{1} \subset \Delta$ denote the finite set of all one-dimensional cones of $\Delta$. We restrict ourselves to toric varieties $X(\Delta)$ which are 
(\romannumeral 1) smooth and (\romannumeral 2) projective; translated to the language of fans these conditions
correspond to the restrictions (\romannumeral 1) $\Delta_{1}$ forms a lattice basis of $N$ in the 
real vector space $  N_{\RR}$ and (\romannumeral 2) the support of the fan 
$|\Delta|$ is all of $N_{\RR}$.  

There is a bijection between cones of the fan $\Delta$ and the orbits of the torus $T$ in $X(\Delta)$: the torus orbit
$O_{\sigma} \subset X(\Delta)$, corresponding to a cone $\sigma \subset \Delta$, is the set of all points in $X(\Delta)$
which are stabilized by the subtorus $\Hom_{gp}(M/M(\sigma), \mb{G}_{m})$ of $T$; where $M(\sigma) \defeq \sigma^{\perp}
\cap M $ is the subspace of characters trivial on the cone $\sigma$. It turns out that the set $O_{\sigma}$ is a Zariski
open subset of its Zariski closure in $X(\Delta)$. This shows that, among other things, that $\dim_\RR (O_{\sigma}) =
\mt{codim}_{\RR}(\sigma)$. In particular cones of maximal dimension in $\Delta$ correspond to $T$-fixed points and the cones of codimension one correspond to $T$- stable curves. 
The identification of orbits and cones imply $O_{\sigma} = T/T_{\sigma}$ and hence
\[
K_{T,*}(O_{\sigma}) = K_{T,*} \otimes R(T/T_{\sigma}) = K_{*}(k) \otimes R(T_{\sigma})= K_{*}(k)[M(\sigma)].
\]

In the light of Theorem \ref{VV}, we conclude that natural restriction map   
\begin{equation}\label{E:KTHTORIC}
K_{T, *}(X(\Delta)) \hookrightarrow \prod_{\sigma \in \Delta_{\max}} 
K_{*}(k) \otimes R(T_{\sigma})
\end{equation}
where $\Delta_{\max}$ is the set of maximal cones in $\Delta$ is an injective map of $R(T)$ modules. Moreover the image
is characterized by the collection of elements 
\[
(a_{\sigma}) \in \prod_{\sigma \in \Delta_{\max}} K_{*}(k) \otimes R(T_{\sigma})
\] 
which satisfy: for any two maximal cones $\sigma_{1}$ and $\sigma_{2}$ of $\Delta$ the restrictions of
$a_{\sigma_{1}}$ and $a_{\sigma_{2}}$ to $R(T_{\sigma_{1}\cap \sigma_{2}})$ coincide.

Concretely, for any maximal cone $\sigma$ of $\Delta$ let us identify $T_{\sigma}$ with $T$, and $T_{\sigma_{1} \cap
  \sigma_{2}}$ with the codimension one subtorus $\Ker(\chi_{\sigma_{1} \cap \sigma_{2}}) \subset T$, where
$\chi_{\sigma_{1} \cap \sigma_{2}}$ is the unique generator of $M_{T_{\sigma_{1} \cap \sigma_{2}}}^{\perp}$. Then an element
$(f_{\sigma})_{\sigma \in \Delta_{\max}}$
in the right hand side of eqn.(\ref{E:KTHTORIC}) belongs to $K_{T, *}(X(\Delta))$ if and only if it satisfies the
condition \begin{equation} \label{E:toriccond}
  f_{\sigma_{1}} - f_{\sigma_{2}} = 0 \mod \; (1- \chi_{\sigma_{1} \cap \sigma_{2}})
\end{equation}
for any pair of intersecting maximal cones $\sigma_{1}$ and $\sigma_{2}$. 

A presentation of $K$-theory in this form will be called \emph{GKM} presentation.

\begin{Remark}\label{R:Another result of VV}
We note in passing that one can consider a sheafified version of equivariant $K$-groups on toric varieties and use the
techniques developed in \cite{Thomason1985} to extend these results to general, not necessarily smooth or projective, toric varieties. 

In particular one shows that for affine (even non-smooth!) toric varieties $U_{\sigma}$ there is a formula
\[K_{T,*}(U_{\sigma}) = R(T_{\sigma})\otimes K_{\ast}(k)\] where $T_{\sigma}$ is the stabilizer of a geometric point.
For a general toric variety $X(\Delta)$ using the natural cover by affine opens sets $\{U_{\sigma}\}_{\sigma \in
  \Delta_{\max}}$ and the exactness of the complex  
\begin{equation}\label{E:trres}
  \begin{tikzcd}
0 \arrow[r] &  K_{T,*}(X) \arrow[r] & \bigoplus_{\sigma \in \Delta_{\max}} K_{T,*}( U_\sigma) \arrow[r, "\partial"] & \bigoplus_{ \sigma,\tau \in \Delta_{\max},\, \sigma \cap \tau \neq \emptyset} K_{T,*}( U_{\sigma \cap \tau}) \\
   {}   &  {}     & (f_{\sigma})_{\sigma \in \Delta_{\max}}  \arrow[r, mapsto, "\partial"] & (f_\tau |_{U_{\sigma \cap \tau}} - f_\sigma |_{U_{\sigma \cap \tau}})_{ \sigma,\tau \in 
     \Delta_{\max},\, \sigma \cap \tau \neq \emptyset}
 \end{tikzcd}
\end{equation}
one can generalize eqn.(\ref{E:toriccond}); see \cite{Au0} for details. 
\end{Remark}

There is another description in-terms of generators and relations for the torus equivariant $K$-theory of smooth projective
toric varieties $X(\Delta)$ called the \emph{multiplicative Reisner-Stanley} (RS) presentation. These presentations were extensively
studied in the context of equivariant cohomology of regular embeddings by Biffet, De Concini and Procesi in \cite{BDP90}.

\begin{Proposition}[The Reisner-Stanley presentation]\label{P:str}
Let $X(\Delta)$ be a smooth toric variety. If $\Delta_1$ denotes the set of all one dimensional 
cones of $\Delta$, then there is an isomorphism of $K_{*}(k)$-algebras
\begin{equation}\label{E:KTHRS}
i: \frac{K_{*}(k)[x_{\rho}^{\pm 1}]}{(\prod_{\rho \in S}(x_{\rho} -1))} \simeq K_{T,*}(X(\Delta)), 
\end{equation}
where the product in the quotient is taken over all subsets $S\subseteq \Delta_{1}$ satisfying the condition: 
\begin{equation}\label{p:cond}
 \mt{the elements of $S$ are not all contained in a maximal cone in $\Delta$.}  
\end{equation}  
\end{Proposition}

\begin{proof}[Sketch of the Proof]
Given any one-dimensional cone $\rho \in \Delta_{1}$, let $v_\rho$  denote the generator of the (rank-one) monoid $\N \cdot \rho
\cap N$.  Let us denote the one dimensional faces of a cone $\sigma$ in the fan $\Delta$ by $\rho_{1,\sigma}, \rho_{2,
  \sigma}, \ldots, \rho_{k, \sigma}$. The fact that $X(\Delta)$ is smooth and projective ensures that for any maximal cone
$\sigma$ the vectors $v_{\rho_{1,\sigma}}, v_{\rho_{2, \sigma}}, \ldots, v_{\rho_{k, \sigma}}$ form an integral basis of
$N$ and the dual generators $ v_{\rho_{i,\sigma}}^{\vee}$ form an integral basis of $M$.

Given any one-dimensional cone $\rho \in \Delta_{1}$ and any $\sigma \in \Delta_{\max}$ consider the assignment
\[
 u_{\rho}^{\sigma} 
 \defeq \begin{cases} 
1 &  \mt{ if } \rho^{\prime} \mt{ is not a face of the cone } \sigma^{\prime},  \\
v_{\rho_{i, \sigma}}^{\vee} & \mt{ if } \rho = \rho_{i}\mt{ is a face of the cone } \sigma 
 \end{cases} 
\]
in $R(T)= \Z[M]$.

The map $i$ in eqn.(\ref{E:KTHRS}) is defined by mapping
\begin{equation}\label{E:RSDEF}
 x_{\rho} \mapsto (u_{\rho}^{\sigma})_{\sigma \in \Delta_{\max}} \in \prod_{\sigma \in \Delta_{\max}} R(T_{\sigma})
\end{equation}
as $\rho$ varies over all one dimensional cones of $\Delta$.

We refer the reader to \cite[Theorem 6.4]{VezzosiVistoli03} which shows that the assignment in eqn.(\ref{E:RSDEF})
defines the correct image. 
\end{proof}

\begin{Remark}
 We point out that the presentation, given by eqn.(\ref{E:KTHRS}), doesn't explicitly show the $R(T)$-module structure on
 the left hand side. It can be recovered from the description of the map $i$ given by (eqn.\ref{E:RSDEF}) above on a
 case by case basis. 
\end{Remark}

Next, we will consider four cases which will be used in the sequel. 
\subsection{Toric $\mathbb{P}^{1}$} \label{SS:p1comp}
The fan of $\mb{P}^{1}$ as a toric variety is shown in Figure \ref{F:P1}. In this case the maximal cells are one-dimensional. 
However, to be consistent with the notation used before, we continue to use $\sigma_{i}$ to denote the maximal cones and
$\rho_{i}$ to denote the one dimensional rays. The underlying torus $T$ in this case is one dimensional so $R(T)=
\Z[\chi^{\pm}]$. We have the formula 
\[ \prod_{\sigma \in \Delta_{\max}} K_{*}(k) \otimes R(T_{\sigma})= K_{*}(k)[\chi^{\pm}]|_{\sigma_{1}} \times K_{*}(k)[\chi^{\pm}]|_{\sigma_{2}}. \]

\begin{figure}[htbp]
\begin{center}
\begin{tikzpicture}[scale = 1.0, transform shape]
\begin{scope}

\shade[fill=gray, draw= gray] (-5,-1)--(5,-1)--(5,1)--(-5,1) --(-5,-1);

\draw (0,0) (0,0.3) node {(0,0)}; 
\filldraw[fill=gray , draw=black ] circle (2 pt);
\path[draw, thick,<->] (3,0) -- (-3,0) ;
\draw (3.8, 0) node {$\rho_{1}= \sigma_{1}$};
\draw (2,-0.3)  node {(1,0)};
\filldraw (2,0) circle (1.5pt);
\draw (-3.8, 0) node {$\rho_{2}= \sigma_{2}$};
\draw (-2,-0.3)  node {(1,0)};
\filldraw (-2,0) circle (1.5pt);

\end{scope}
\end{tikzpicture}
\end{center}
\caption{The fan of $\mb{P}^{1}$.}
\label{F:P1}
\end{figure}

\subsubsection*{GKM Presentation} We note that the intersection of the maximal cones $\sigma_{1} \cap \sigma_{2}$ is the zero cone whose stabilizer is the trivial group. 
As a result, the GKM description of the $K$ groups is 
\begin{multline*}
K_{T,*}(\mb{P}^{1}) = \left \{ (f_{1}, f_{2}) \in K_{*}(k)[\chi^{\pm}] \times K_{*}(k)[\chi^{\pm}] \; | \; \text{the constant term of}\; f_{1} \right. \\ 
\left . = \text{the constant term  of} \; f_{2} \right \}.
\end{multline*}

\subsubsection*{RS Presentation}  
The Reisner-Stanley presentation is given by \[ \frac{K_{*}(k) [x_{\rho_{1}}^{\pm}, x_{\rho_{2}}^{\pm}]}{ (x_{\rho_{1}}-1)(x_{\rho_{2}}-1) }.\]
Using Proposition \ref{P:str} we see that the generators are mapped to 
\[ x_{\rho_{1}} \mapsto (\chi, 1)  \text{ and } x_{\rho_{2}} \mapsto (1,\chi).\]
As a $K_{*}(k)[\chi^{\pm}]$ module the action of $\chi$ on the generators is multiplication by $x_{\rho_{1}}x_{\rho_{2}}$.   
 
\subsection{Surfaces} \label{SS:toric comp}
We consider the case of projective plane and the Hirzebruch surfaces (including $\PP^{1} \times \PP^{1}$). The
underlying torus $T$ is two dimensional and in co-ordinates  $R(T) = \Z[ \chi_{1}^{\pm}, \chi_{2}^{\pm}]$. We identify
$K_{*}(k) \otimes R(T)$ with $K_{*}(k)[\chi_{1}^{\pm} , \chi_{2}^{\pm}]$. We use  $K_{*}(k)[\chi^{\pm}]$ to denote
$K_{*}(k)[\chi_{1}^{\pm}, \chi_{2}^{\pm}]$, and with this notation, the right hand-side of eqn.(\ref{E:KTHTORIC}) is  
\[
\prod_{\sigma_{ij} \in \Delta_{\max}} K_{*}(k) \otimes R(T_{\sigma_{ij}}) = K_{*}(k)[\chi^{\pm}]|_{\sigma_{12}}  \times \ldots \times  K_{*}(k)[\chi^{\pm}]|_{\sigma_{ij}}
\]
with a diagonal $K_{*}(k)[\chi^{\pm}]$ action.
When no confusion is likely, we denote the ring $\prod_{\sigma_{ij} \in \Delta_{\max}} K_{*}(k) \otimes
R(T_{\sigma_{ij}}) $ will be denoted by $K_{*}(k)[\chi^{\pm}_{\sigma}]$. 

\subsubsection{The projective plane: $\mathbb{P}^{2}$} \label{sb:p2comp}

\begin{figure}
\begin{minipage}{0.4\textwidth}

\begin{center}
\begin{tikzpicture}[scale = 0.8, transform shape]
\begin{scope}
\shade[fill=gray, draw=gray] (0,0)--(3,0)--(3,3)--(0,3)--(-3,2.8)--(-2.0,-2.0)--(0,-3)--(3,-3)--(3,0);
\path[draw, thick, ->] (0,0) -- (0,3) ;
\path[draw, thick,->] (0,0) -- (3,0) ;
\path[draw, thick,->] (0,0) -- (-2,-2) ;

\draw (3.3, 0) node {$\rho_{1}$};
\draw (1,-0.3)  node {(1,0)};
\filldraw (1,0) circle (1.5pt);
\draw (2, 2) node {$\sigma_{12}$};
\draw (0, 3.3) node {$\rho_{2}$};
\draw (-0.5,1.0)  node {(0,1)};
\filldraw (0,1) circle (1.5pt);
\draw (-2, 2) node {$\sigma_{23}$};
\draw (-2.3, -2.3) node {$\rho_{3}$};
\draw (-0.5,-1.4)  node {(-1,-1)};
\filldraw (-1,-1) circle (1.5pt);
\draw (2,-2) node {$\sigma_{13}$};

\end{scope}
\end{tikzpicture}
\caption{The fan of $\mb{P}^{2}$.}
\label{F:P2}
\end{center}
\end{minipage}
\hfill
\begin{minipage}{0.4\textwidth}
\begin{center}
\begin{tikzpicture}[scale = 0.7, transform shape]
\begin{scope}
\shade[fill= gray, draw = gray] (0,0)--(3,0)--(3,3)--(0,3)--(-3, 3)--(-3,0)--(-3,-3)--(3,-3)--(3,0);

\path[draw, thick, <->] (0,3) -- (0,-3) ;
\path[draw, thick,<->] (3,0) -- (-3,0) ;
\draw (3.3, 0) node {$\rho_{1}$};
\draw (1,-0.3)  node {(1,0)};
\filldraw (1,0) circle (1.5pt);
\draw (2, 2) node {$\sigma_{12}$};
\draw (0, 3.3) node {$\rho_{2}$};
\draw (-0.4,1.0)  node {(0,1)};
\filldraw (0,1) circle (1.5pt);
\draw (-2, 2) node {$\sigma_{23}$};
\draw (-3.3, 0) node {$\rho_{3}$};
\draw (-1,-0.3)  node {(1,0)};
\filldraw (-1,0) circle (1.5pt);
\draw (0,-3.3) node {$\rho_{4}$};
\draw (-0.4,-1.0)  node {(0,1)};
\filldraw (0,-1) circle (1.5pt);
\draw (-2, -2) node {$\sigma_{34}$};
\draw (2,-2) node {$\sigma_{14}$};

\end{scope}
\end{tikzpicture}
\caption{The fan of $\mb{P}^{1} \times \mb{P}^{1}$.}
\label{F:P1P1}
\end{center}
\end{minipage}
\end{figure}

The fan of $\mb{P}^{2}$ is shown in Figure \ref{F:P2}. 

\subsubsection*{GKM presentation}
Explicitly the left-hand-side of eqn.(\ref{E:KTHTORIC}) is given by
\begin{align*}
       K_{*}(k) \otimes R(T_{\rho_{1}})& = K_{*}(k)[\chi_{1}^{\pm}]\\
         K_{*}(k) \otimes R(T_{\rho_{2}}) & = K_{*}(k)[\chi_{2}^{\pm}]\\
         K_{*}(k) \otimes R(T_{\rho_{3}}) & = K_{*}(k)[(\chi_{1}\chi_{2})^{\pm}].
        \end{align*}   

An element $(f_{1}, f_{2}, f_{3})\in K_{*}(k)[\chi^{\pm}_{\sigma}]$ belongs to $K_{T,*}(\mb{P}^{2})$ if and only if $f_{1} - f_{2} =0 \mod (1-\chi_{1})$, $f_{1} - f_{3} = 0 \mod (1-\chi_{2})$ and $f_{2} - f_{3} = 0 \mod (1- \chi_{1}^{-1}\chi_{2})$.

\subsubsection*{RS presentation}
The Reisner-Stanley presentation is given by \[ \frac{K_{*}(k) [x_{\rho_{1}}^{\pm}, x_{\rho_{2}}^{\pm}, x_{\rho_{3}}^{\pm}]}{ (x_{\rho_{1}}-1)(x_{\rho_{2}}-1)(x_{\rho_{3}}-1) }\]
where the generators $x_{\rho_{i}}$ are mapped, via eqn.(\ref{E:KTHRS}), to    
\begin{align*}
        x_{\rho_{1}} & \mapsto (\chi_{1},\; 1,\; \chi_{1}\chi_{2}^{-1})\\
        x_{\rho_{2}} & \mapsto (\chi_{2},\; \chi_{1}^{-1}\chi_{2}, \;1)\\
        x_{\rho_{3}} & \mapsto (1, \;\chi_{1}^{-1}, \;\chi_{2}^{-1}) 
        \end{align*}   
As a $K_{*}(k)[\chi^{\pm}]$ module, $\chi_{1}$ action is multiplication by $x_{\rho_{1}}x_{\rho_{3}}^{-1}$ and $\chi_{2}$ action is multiplication by $x_{\rho_{2}}x_{\rho_{3}}^{-1}$.

\subsubsection{The surface: $\mb{P}^{1} \times \mb{P}^{1}$} \label{sb:p1p1}

The fan of $\mb{P}^{1}\times \mb{P}^{1}$ is shown in Figure \ref{F:P1P1}. 
\subsubsection*{GKM presentation}
The left-hand-side of eqn.(\ref{E:KTHTORIC}) is given by 
\begin{align*}
       K_{*}(k) \otimes R(T_{\rho_{1}})& = K_{*}(k)[\chi_{1}^{\pm}]\\
         K_{*}(k) \otimes R(T_{\rho_{2}}) & = K_{*}(k)[\chi_{2}^{\pm}]\\
         K_{*}(k) \otimes R(T_{\rho_{3}}) & = K_{*}(k)[\chi_{1}^{\pm}]\\
          K_{*}(k) \otimes R(T_{\rho_{4}}) & = K_{*}(k)[\chi_{2}^{\pm}].
        \end{align*}   
The relations are given by: an element  $(f_{1}, f_{2}, f_{3}, f_{4})\in K_{*}(k)[\chi^{\pm}_{\sigma}]$ belongs to $K_{T,*}(\mb{P}^{1}\times \mb{P}^{1})$ if and only if $f_{1} - f_{2} = 0 \mod (1-\chi_{1})$, $f_{2} - f_{3} = 0 \mod (1-\chi_{2})$, $f_{3} - f_{4} = 0 \mod (1-\chi_{1})$ and $f_{4} - f_{1} = 0 \mod (1-\chi_{2})$.

\subsubsection*{RS presentation}The Reisner-Stanley presentation is given by the formula \[ \frac{K_{*}(k) [x_{\rho_{1}}^{\pm}, x_{\rho_{2}}^{\pm}, x_{\rho_{3}}^{\pm}, x_{\rho_{4}}^{\pm}]}{ ((x_{\rho_{1}}-1)(x_{\rho_{3}}-1), (x_{\rho_{2}}-1)(x_{\rho_{4}}-1))}.\]
 where the generators $x_{\rho_{i}}$ are mapped as follows 
\begin{align*}
        x_{\rho_{1}} & \mapsto (\chi_{1},\;1, \;1,\; \chi_{1})\\
        x_{\rho_{2}} & \mapsto (\chi_{2},\; \chi_{2},\; 1, \;1)\\
        x_{\rho_{3}} & \mapsto (1,\; \chi_{1}^{-1}, \; \chi_{1}^{-1}, \;1) \\
        x_{\rho_{4}} & \mapsto (1,\; 1,\; \chi_{2}^{-1},\; \chi_{2}^{-1}).
        \end{align*}   
In this case, as a $K_{*}(k)[\chi^{\pm}]$ module, $\chi_{1}$ action is multiplication by $x_{\rho_{1}}x_{\rho_{3}}^{-1}$ and $\chi_{2}$ action is multiplication by $x_{\rho_{2}}x_{\rho_{4}}^{-1}$.

\subsubsection{The Hirzebruch surfaces: $\mb{F}_{n}$, $n > 1$} 
\label{sb:fn}
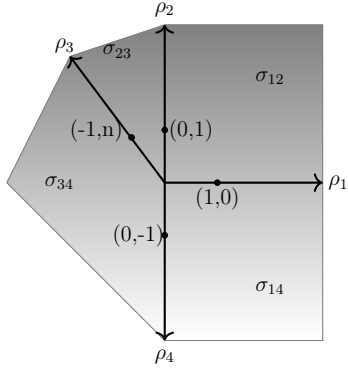
\begin{figure}[htbp] 
\begin{center}
\begin{tikzpicture}[scale = 0.7, transform shape]
\begin{scope}
\shade[fill=gray, draw=gray] (0,0)--(3,0)--(3,3)--(0,3)--(-1.8, 2.4)--(-3,0)--(0,-3)--(3,-3)--(3,0);
\path[draw, thick, ->] (0,0) -- (0,3) ;
\path[draw, thick,->] (0,0) -- (3,0) ;
\path[draw, thick,->] (0,0) -- (-1.8,2.4) ;
\path[draw, thick,->] (0,0) -- (0,-3) ;
\draw (3.3, 0) node {$\rho_{1}$};
\draw (1,-0.3)  node {(1,0)};
\filldraw (1,0) circle (1.5pt);
\draw (2, 2) node {$\sigma_{12}$};
\draw (0, 3.3) node {$\rho_{2}$};
\draw (0.5,1.0)  node {(0,1)};
\filldraw (0,1) circle (1.5pt);
\draw (-0.9, 2.5) node {$\sigma_{23}$};
\draw (-1.9, 2.6) node {$\rho_{3}$};

\draw (-1.3,1.0)  node {(-1,n)};
\filldraw (-0.63, .86) circle (1.5pt);

\draw (0,-3.3) node {$\rho_{4}$};
\draw (-0.5,-1.0)  node {(0,-1)};
\filldraw (0,-1) circle (1.5pt);
\draw (-2, 0) node {$\sigma_{34}$};
\draw (2,-2) node {$\sigma_{14}$};

\end{scope}
\end{tikzpicture}
\caption{ The fan of Hirzebruch surface $\mb{F}_{n}$.}
 \label{F:Fn}
\end{center}
\end{figure}

The fan of a Hirzebruch surface $\mb{F}_{n}$ is shown in Figure \ref{F:Fn}.
\subsubsection*{GKM presentation}
Using the co-ordinates $\chi_{1}$ and $\chi_{2}$ on the torus the left-hand-side of eqn.(\ref{E:KTHTORIC}) is given by
\begin{align*}
       K_{*}(k) \otimes R(T_{\rho_{1}})& = K_{*}(k)[\chi_{1}^{\pm}]\\
         K_{*}(k) \otimes R(T_{\rho_{2}}) & = K_{*}(k)[\chi_{2}^{\pm}]\\
         K_{*}(k) \otimes R(T_{\rho_{3}}) & = K_{*}(k)[(\chi_{1}^{-1}\chi_{2}^{n})^{\pm}]\\
          K_{*}(k) \otimes R(T_{\rho_{4}}) & = K_{*}(k)[\chi_{2}^{\pm}].
        \end{align*}   
An element $(f_{1}, f_{2}, f_{3}, f_{4})\in K_{*}(k)[\chi^{\pm}_{\sigma}]$ belongs to $K_{T,*}(\mb{F}_{n})$ if and only if $f_{1} - f_{2} = 0 \mod (1-\chi_{1})$, $f_{2} - f_{3} = 0 \mod (1-\chi_{1}^{n}\chi_{2})$, $f_{3} - f_{4} = 0 \mod (1-\chi_{1})$ and $f_{4} - f_{1} = 0 \mod (1-\chi_{2})$. 

\subsubsection*{RS presentation}
The Reisner-Stanley presentation is given by \[ \frac{K_{*}(k) [x_{\rho_{1}}^{\pm}, x_{\rho_{2}}^{\pm}, x_{\rho_{3}}^{\pm}, x_{\rho_{4}}^{\pm}]}{ ((x_{\rho_{1}}-1)(x_{\rho_{3}}-1), (x_{\rho_{2}}-1)(x_{\rho_{4}}-1))}.\]

The generators $x_{\rho_{i}}$ are mapped as follows
\begin{align*}
        x_{\rho_{1}} & \mapsto (\chi_{1},\;1,\; 1,\; \chi_{1})\\
        x_{\rho_{2}} & \mapsto (\chi_{2},\;\chi_{1}^{n}\chi_{2},\; 1,\; 1)\\
        x_{\rho_{3}} & \mapsto (1,\; \chi_{1}^{-1}, \;\chi_{1}^{-1},\; 1) \\
        x_{\rho_{4}} & \mapsto (1,\; 1,\; \chi_{1}^{-n}\chi_{2}^{-1},\; \chi_{2}^{-1}).
        \end{align*}   
As a $K_{*}(k)[\chi^{\pm}]$ module, $\chi_{1}$ action is multiplication of $x_{\rho_{1}}x_{\rho_{3}}^{-1}$ and $\chi_{2}$ action is multiplication by $x_{\rho_{2}}x_{\rho_{3}}^{n}x_{\rho_{4}}^{-1}$.

%
%
%
%

\section{Torus actions on Spherical varieties}\label{S:Preliminaries}
In this section we will analyze the irreducible components of the \emph{fixed point locus} $X^{S}$ of diagonalizable
subgroups $S$ of a reductive group $G$ and a $G$-variety $X$. 

More precisely, for any $G$-variety $X$ and any closed subgroup $S \subset G$ we consider the functor $\underline{X^{S}}$ which associates to any affine $k$-scheme $A$ the set \[ \underline{X^{S}}(A)= \left \lbrace x \in X(A) \; | s\cdot x = x \; \text{ for any } s \in S(A) \right \rbrace.\]
The functor $\underline{X^{S}}$ is a representable closed sub-functor of $\underline{X}$. The \emph{fixed point locus}
of $S$ on $X$ is the closed subscheme $X^{S}$ of $X$, representing $\underline{X}^{S}$, with the reduced scheme
structure.\footnote{It was pointed out by Brion that a result of Fogarty implies that $X^{S}$ is already smooth; see \cite{fogarty}.}  

We recall some standard results about reductive algebraic groups over fields. The main reference for these results is \cite{Borel}.
 \begin{Definition}
Let $S$ be any torus in $G$. Then $S$ is called \emph{regular} if $S$ contains a regular element\footnote{Recall an
  element $g \in G$ is regular if the dimension of the centralizer $C_{G}(g)$ is minimal.} of $G$ and $S$ is called \emph{singular} if $S$ is contained in infinitely many Borel subgroups of $G$.
\end{Definition}

\begin{Remark} 
We note that Borel, in \cite[\Rnumeral{4}, \S 13.7, Cor. 2]{Borel}, considers a more general notion of a \emph{semi-regular} torus but it turns out that for reductive
groups a semi-regular torus is regular. So, for reductive groups, we have a dichotomy: a torus $S$ is either regular or singular. 
\end{Remark}

A maximal torus is always regular. So a singular torus $S$, contained in a maximal torus $T$, is a subtorus of
codimension at-least one. All codimension one singular tori, contained in a fixed maximal torus $T$, correspond to the
neutral component of $\Ker(\alpha)$ where $\alpha$ is a root of $G$ (with respect to $T$).    

\subsubsection*{Centralizers} Suppose $S$ is a codimension one subtorus of a maximal torus $T \subset G$. Then the
centralizer $C_{G}(S)$ is a connected reductive group. If $S$ is regular, then $C_{G}(S) = T$, and if $S$ is singular,
then the semisimple part $C_{G}(S)^{ss}$ has rank one. If $B$ is any Borel subgroup of $G$
containing $S$ then $C_{B}(S)$ maps onto a Borel subgroup of $C_{G}(S)^{ss}$ and conversely any Borel subgroup of
$C_{G}(S)^{ss}$, containing the image of $S$, is the image of a group of the form $C_{B}(S)$.   

\subsection{Spherical Varieties} \label{S:fixed point loci}
\begin{DefProp}[See \cite{Brion86}]
An irreducible $G$-variety $X$ is called a \emph{spherical variety} if any of the following equivalent conditions hold.
\begin{enumerate}
\item For any Borel subgroup $B$ of $G$, the only $B$-invariant rational functions on $X$ are constant functions. 
\item The minimal codimension of a $B$-orbit in $X$ is zero i.e. $X$ has a dense open $B$-orbit. 
\item $X$ has finitely many $B$-orbits. 
\end{enumerate}
A \emph{homogeneous spherical variety} is a homogeneous space $G/H$ which is also a spherical $G$-variety.  
\end{DefProp}

The open $G$-orbit of a spherical variety $X$ is a homogeneous spherical variety $G/H$ and this is equivalent to the
condition that the subvariety $B \cdot H$ is open in $G$. The subgroups $H \subset G$ satisfying this property are
called \emph{spherical subgroups}. It follows from Condition (3) and some additional arguments\footnote{It is not
  immediate that the closure is a normal variety.} that the closure of a $G$-orbit inside a spherical variety is a spherical variety.  

The set of $B$-eigenvalues of the eigenvectors in $k(X)$, the field of rational functions on $X$, is a sublattice of the
lattice of characters of $B$. The rank of this lattice is a fundamental invariant of the spherical variety $X$ and we will denote it by $r(X)$.   

\subsubsection{Torus action on spherical varieties} Consider a spherical $G$-variety $X$. Throughout this section we fix
a Borel subgroup $B$, a maximal torus $T \subset B$ and let $S$ denote any codimension one subtorus of $T$. 

\begin{Lemma}
The set of $T$-fixed points of a spherical variety $X$ is finite. 
\end{Lemma}

\begin{proof}
The spherical variety $X$ decomposes as a finite union of $G$-orbits. So the Lemma is equivalent to showing that: if $x \in X$ is
any torus fixed point then the $G$-orbit $G \cdot x$ has finitely many torus fixed points.

The latter statement follows from Lemma 2.2 of \cite{decsp}, which shows that a $G$-homogeneous space
$G \cdot x$ has finitely many $T$ fixed points if the stabilizer of $x$ contains $T$.

\end{proof}

\begin{Remark}
If $X$ is complete then Borel Fixed Point Theorem shows that $X^{T}$ is nonempty. In particular, $X^T$ has at least $\dim(X) +1$ points, \cite[Theorem 10.2, \Rnumeral{4}]{Borel}. 
\end{Remark}

Now consider the fixed point locus $X^{S}$. It is stable under the action of the centralizer $C_G(S)$ and the structure of $X^{S}$, somewhat unsurprisingly, depends on whether $S$ is regular or singular.

\begin{Lemma}\label{L:single B orbit}
Let $x$ be a $S$-fixed point of $X$.
The intersection of $B\cdot x$ with the fixed point loci $X^S$  
is the $C_{B}(S)$ orbit $C_B(S) \cdot x$.
\end{Lemma}

\begin{proof}
It is clear that $C_B(S) \cdot x \subseteq B\cdot x \cap X^S$. 

Let us fix a realization, \emph{\`{a} la Springer} \cite[Chapter 8]{SpringerBook}, of $G$ with respect to the root system $\Phi \defeq \Phi(G,T)$ and $\Phi^{+}$ is the subset of positive roots. This gives us a family of homomorphisms $\varphi_{\alpha} : \mb{G}_{a} \rightarrow G$, indexed by $\alpha \in \Phi$, such that for any $t \in T$ we get $t \cdot \varphi_{\alpha} \cdot t^{-1} = u_{\alpha}(\alpha(t)\cdot x)$.  The Borel subgroup $B$ admits a decomposition $B=TU$ and $U$ is generated by the images of $(\varphi_{\alpha})_{\alpha \in \Phi^{+}}$. 
To prove the converse assertion let $y\in B\cdot x \cap X^S$. Then $y= b_{0} \cdot x$ for some $b_{0}\in B$. We write $b_{0} = t_{0} \cdot u_{0}$ where $t_{0} \in T$, $u_{0} \in U$ and we can moreover assume $u_{0} = \varphi_{\alpha}(z_{0})$ for some $z_{0} \in \mb{G}_{a}$. 

We have \[ s \cdot b_{0} \cdot s^{-1} = t_{0} \cdot \varphi_{\alpha}(\alpha(s)z_{0}). \]

When $S$ is a singular torus and $\alpha$ is a root such that $S \subset \mt{Ker}(\alpha)$ then clearly $b_{0} \in
C_{B}(S)$. When $\alpha$ is a nontrivial character of $S$ we have, for any $s \in S$, \[ y = s \cdot y = t_{0} \varphi_{\alpha}(\alpha(s)z_{0}) \cdot x = \lim_{s \rightarrow 0} t_{0} \varphi_{\alpha}(\alpha(s)z_{0}) \cdot x = t_{0} \cdot x. \]

This shows that we can write $y = t_{0} \cdot x$ for some $t_{0} \in T$ and the conclusion follows. 
\end{proof}

\begin{Corollary}\label{C:X^S is spherical}
We continue to use the notation of Lemma \ref{L:single B orbit}.
In this case, if $Y$ is any irreducible component of $X^{S}$ then $Y$ is a spherical $C_G(S)$-variety.
\end{Corollary}
\begin{proof}
The normality of $Y$ follows from \cite{fogarty}. The group $C_{G}(S)$ is connected, so it
stabilizes $Y$.
The ambient variety $X$ is spherical hence it has finitely many $B$-orbits; therefore using the previous lemma we see that $Y$ has finitely many $C_{B}(S)$ orbits. Hence $Y$ is a spherical $C_{G}(S)$ variety. 
\end{proof}

\begin{Corollary}\label{C:jh}
We continue with the notation of Corollary \ref{C:X^S is spherical}. The variety $Y$ has a dimension at most two and as a $C_{G}(S)$-spherical variety it has rank $r(Y)$ at most one.    
\end{Corollary}
\begin{proof}
The dimension of $Y$ is bounded by the dimension of $C_{G}(S)$.
The assertion is clear when the quotient $C_{G}(S)/S$ has dimension one. 
Consider the case when $S$ is singular and hence dim($C_{G}(S)/S$) is greater than one. In this case the
dimension and rank of $Y$ are determined by the dimension and rank of the open $C_{G}(S)/S$-orbit. 

The center $Z_{S}$ of $C_{G}(S)$ is contained in $S$. As a result,                                          
 the action of $C_{G}(S)$ on $Y$ factors through the semisimple part $\mc{G}_{S} \defeq C_{G}(S)^{ss}$. The group
$C_{G}(S)^{ss}$ has semisimple rank one and hence $ \mc{G}_{S}= SL_{2} \; \text{ or } \; PSL_{2}$. Let $\mc{B}_{S}$ denote the
image of $C_{B}(S)$  in $\mc{B}_{S}$ . 

The dimension of $Y$ is bounded by the dimension of $\mc{B}_{S}$ which is two, and the rank of $Y$ is 
bounded by the dimension of the character lattice of $\mc{B}_{S}$ which is at most one. 
This proves the assertion.
\end{proof}

\subsection{Rank one spherical varieties}
A complete classification of rank one spherical varieties is known, see \cite{Ahiezer, Brion1987} for details. In this
section, we drop
the sub-script $S$ from the semi-simple group $\mc{G}_{S}$, the stabilizer $\mc{H}_{S}$ etc. When $\mc{G}$ is $SL_{2}$ or
$PSL_{2}$\footnote{We use $\mc{G}$ to denote the semi-simple rank one group case. The is to distinguish from the arbitrary reductive
  group case, denoted by $G$, in the next section.} the classification is determined by the
spherical subgroups of $\mc{G}$. These rank one compactifications will be used in the next section. We recall the spherical subgroups, equivariant compactifications and the boundary of the open $\mc{G}$ orbit in Table \ref{Table:allcases}.     

\begin{table}[htp]
\rowcolors{2}{gry}{white}
\begin{center}

\begin{tabulary}{0.9\textwidth}{|C|C|C|}

\hline 
Spherical subgroup $\mc{H} \subset \mc{G}$ & Equivariant 
 of $\mc{G}/\mc{H}$  
& Boundary of open $\mc{G}$ orbit
\\
\hline
$B$ = Borel subgroup & $\mb{P}^{1}$ & \\
 
 $T=$ Maximal Torus & $\mb{P}^{1} \times \mb{P}^{1}$ & diagonal $\mb{P}^{1}$ 
 \\
  $N_{G}(T)=$ Normalizer of $T$ & $\mb{P}(\mf{sl}_{2})$ &  conic of nilpotent matrices
 \\
 $C_{n}\ltimes U$,  where $C_{n} = \text{diag}(\zeta, \zeta^{-1})$ for  $\zeta^{n} = 1$ and $U$ is the unipotent subgroup &  $\mb{F}_n = \mb{P}(\mc{O}_{\mb{P}^1} \oplus \mc{O}_{\mb{P}^1}(n))$  & $\mb{P}(\mc{O}_{\mb{P}^1} \oplus 0)$  and 
 $\mb{P}(0 \oplus \mc{O}_{\mb{P}^{1}}(n))$ 
 \\
\hline
\end{tabulary}
\vspace*{1.5ex}
\caption{Equivariant embeddings of semisimple rank one groups.}
\label{Table:allcases}
\end{center}
\end{table}

Let us continue with the notation of the proof of Corollary \ref{C:jh}. Let $y_{0}$ be a generic point in the open $\mc{G}$ orbit and consider the orbit map $\varphi: \mc{G}/\mc{H} \rightarrow Y$ given by $g \mapsto g\cdot y$.  Let $\mc{X}$ denote one of the equivariant compactifications in the Table \ref{Table:allcases} (depending on $\mc{H}$). We consider the birational $\mc{G}$-equivariant map $\varphi: \mc{X} \dashrightarrow Y$. The map $\varphi$ is possibly undefined in a codimension two locus by Zariski's Main Theorem. When $\mc{X}$ is two-dimensional the boundary is one dimensional so the birational map $\varphi$ extends by $\mc{G}$-equivariance to the boundary. 

When $\mc{X}$ is not the surface $\mb{P}^{1} \times \mb{P}^{1}$ the boundary curves are not contractible $(-1)$ curves and $\mc{X}$ is a minimal model. So the map $\varphi$ is necessarily an isomorphism. In the case $\mc{X} = \mb{P}^{1} \times \mb{P}^{1}$ any contractible curve necessarily intersects the dense $\mc{G}$ orbit. So the map $\varphi$ is an isomorphism. 

We get a complete characterization of the geometry of the irreducible two dimensional components of $X^{S}$. Summarizing, we have the following proposition.

\begin{Proposition} \label{P:Brion-characterization}
Suppose $S\subset T$ is a codimension one subtorus, 
and let $Y\subseteq X^S$ be an irreducible. Then up-to isomorphism $Y$ is one of the following varieties. 
\begin{enumerate}[(i)] 
\item a point,
\item a smooth $\mb{P}{^1}$, identified with the complete flag variety $\mc{G}/\mc{B}$.
\item a projective plane on which $C_G(S)$ acts through the projectivization of 
the adjoint action on the Lie-algebra of {\em $SL_2$},  
\item a Hirzebruch surface $\mb{F}_n=\mb{P}( \mc{O}_{\mb{P}^1} \oplus \mc{O}_{\mb{P}^1}( n ))$ for $n > 1$ (or
  $\mb{P}^{1} \times \mb{P}^{1}$) identified with the projectivization of the rank two equivariant vector bundle $\mc{G} \times_{\mc{B}} V $
over $\mc{G}/\mc{B} \rightarrow B $. Here $V = k\cdot \chi \oplus k \cdot \chi^{n} $, for $n \geq 1$, is a two-dimensional $B$-representation
(extended trivially from $T$) and $\chi$ is a generator of the lattice of characters of $T$.
\end{enumerate}
\end{Proposition}

%
%
%
%

\section{$K$-Theory of Spherical Varieties}\label{S:Kthsph}
In this section our goal is to prove Theorem \ref{T:mt1}. This will be achieved in two steps. In the first step, we will
provide a direct proof for the two-dimensional rank-one spherical varieties which appear in Proposition
\ref{P:Brion-characterization}. Then using these results, combined with the work in previous sections, we will prove
the general case.

\subsection{Rank one case}\label{S:Rank one case}
Let us first outline our strategy. Suppose $\mc{X}$ denote any of the two dimensional compactifications listed in
Proposition  \ref{P:Brion-characterization} then it follows that $\mc{X}$ is also a toric variety compactifying a
two-dimensional torus $T$. Let $\mc{T}$ denote a fixed maximal torus of $\mc{G}$ and $\mc{T}$-acts on $\mc{X}$ via the $\mc{G}$-action. 
 
We will show that there is a closed embedding $\iota :\mc{T} \hookrightarrow T$ which makes the Diagram (\ref{E:embed}),
where the vertical maps are the action maps,  commutative.  
\begin{equation}\label{E:embed}
\begin{tikzcd}
  \mc{X} \times \mc{T} \arrow{r}{id \times i} \arrow{d} &  \mc{X} \times T
  \arrow{d}\\
\mc{X} \arrow{r}[equal]{id} & \mc{X}
\end{tikzcd}
\end{equation}
When we have such a map $\iota$, Proposition \ref{L:Edidin's tip} shows that
\begin{equation}\label{E:fds} 
K_{T,*}(\mc{X}) \otimes_{R(T)} R(\mc{T}) \simeq K_{\mc{T},*}(\mc{X}).
\end{equation}
The structure of $K_{T, *}(\mc{X})$ is explicit from the toric computations in Section \ref{SS:toric} and the $R(T)$-module
structure on $R(\mc{T})$ is clear (it depends on $\iota$).

We fix co-ordinates which identifies $\mc{T}$ with the maximal torus 
\[
T= \left\{ 
\begin{pmatrix}
t & 0 \\
0 & t^{-1}
\end{pmatrix}:\ t\in \G_m \right\}.
\]
in $\mc{G}$ and let $\chi$ denote the character \[
\chi \begin{pmatrix} t & 0 \\ 0 & t^{-1}\end{pmatrix} = t\] 
which is a generator of the character lattice of $\mc{T}$ and hence $R(\mc{T}) = \Z[\chi^{\pm}]$. Recall that, using
co-ordinates $\chi_{1}, \chi_{2}$, we have $R(T) = \Z[\chi_{1}^{\pm}, \chi_{2}^{\pm}]$. 

The key to constructing the embedding $\iota: \mc{T} \hookrightarrow T$ is the fact that for all varieties, listed in
Proposition \ref{P:Brion-characterization}, the complement of the open $\mc{G}$- orbit is a $\mc{G}$-stable hence
$\mc{T}$-stable hypersurface. This imposes enough restrictions on $\iota$ to recover the map.

\subsubsection{The case of $\mc{X} = \mb{P}(\mf{sl}_2)$.}\label{SS:p2} 
We identify the coordinates $[x_{0}: x_{1}: x_{2}]$ 
on $\mb{P}({\mf{sl}_2})$ with trace zero matrices: 
$\begin{pmatrix} x_{0} & x_{1} \\ x_{2} & -x_{0} \end{pmatrix}$ 
modulo scalars.
The torus $T$ acts on $\mb{P}^{2}$ is by $t \cdot[x_{0}: x_{1}: x_{2}] = [x_{0}: \chi_{1}(t)\cdot x_{1} : \chi_{2}(t) \cdot x_{2}]$ and the torus $\mc{T}$ 
acts on $\mf{sl}_{2}$ by conjugation action.
The complement of the open $\mc{G}$ orbit is isomorphic 
to the quadric $x_{0}^{2} + x_{1}x_{2} = 0$.
Hence the embedding $\iota: \mc{T} \hookrightarrow T$ is given by $t \mapsto (t, t^{-1})$.
The induced map on the representation rings is given by 
\begin{align*}
 \Z[\chi_{1}^{\pm}, \chi_{2}^{\pm}] \rightarrow & \Z[\chi^{\pm}] \\
   \chi_{1} \mapsto & \chi \\
   \chi_{2} \mapsto & \chi^{-1}.
\end{align*}
The GKM presentation follows immediately from Proposition \ref{L:Edidin's tip}.

\begin{Proposition}\label{P:P2K}
The $\mc{T}$ equivariant $K$-theory of $\mb{P}^{2}$ is given by $3$-tuples $(f_{1}, f_{2}, f_{3})$  where each $f_{i}$ (for $i = 1,2,3$) belongs to the ring $K_{*}(k)[ \chi^{\pm}]$, which satisfies the relations
\[
f_{1} -f_{2} = 0 \mod (1-\chi), \; f_{1} -f_{3} = 0 \mod (1-\chi), \;\text { and } f_{2} -f_{3} = 0 \mod (1-\chi^{2}).
\]
\end{Proposition}

\subsubsection{The case of $\mc{X} = \mb{P}^{1} \times \mb{P}^{1}$.} \label{SS:p1p1}
The diagonal $\mb{P}^{1}$ in $\mc{X}$ is the boundary of the open 
$\mc{G}$-orbit. 
The diagonal intersects the dense $T$ orbit. 
As a result, the only map that preserves the boundary of the compactification 
is the diagonal map. The induced map on the representation rings is given by 
\begin{align*}
 \Z[\chi_{1}^{\pm}, \chi_{2}^{\pm}] \rightarrow & \Z[\chi^{\pm}] \\
   \chi_{1} \mapsto & \chi \\
   \chi_{2} \mapsto & \chi.
\end{align*}
The GKM presentation follows immediately from Proposition \ref{L:Edidin's tip}.

\begin{Proposition}\label{P:P1P1K}
The $\mc{T}$ equivariant $K$-theory of $\mb{P}^{1} \times \mb{P}^{1}$ is given by a $4$-tuple $(f_{1}, f_{2}, f_{3}, f_{4})$,
where each $f_{i}$ (for $i = 1,2,3,4$) belongs to the ring $K_{*}(k)[ \chi^{\pm}]$. Moreover all such tuples satisfy the
relation
\[
f_{i} -f_{j} = 0 \mod (1- \chi) \; \text { for } i,j  \in \{ 1,2,3,4\}.
\]
\end{Proposition}

\subsubsection{The case of $\mc{X} = \mb{F}_{n}$.} \label{SS:fn}
Let $\mc{B}^{-}\subset \mc{G}$ 
denote the Borel subgroup of lower triangular matrices. 
Consider the character $\phi_n$ of $\mc{B}^-$ defined by  
\[
\phi_{n} 
\begin{pmatrix} a & 0 \\ b & a^{-1} \end{pmatrix} =
a^n;
\]
so $\phi_{n} = \chi^{n}$.
The group $\mc{B}^{-}$ acts on $\mb{P}^{1}$ by  $b \cdot [x_{0}: x_{1}] = [x_{0}: \phi_{n}(b) x_{1}]$ and the associated
bundle $\mc{G} \times _{\mc{B}^{-}} \mb{P}^{1}$ , with the canonical fibration to $\mc{G}/\mc{B}^{-} = \PP^{1}$ is identified with the Hirzebruch surface $\mb{F}_{n}$. 

The map $\iota: \mc{T} \rightarrow T$ is given by $t \mapsto (t, t^{n})$.

The induced map on the representation rings is given by 
\begin{align*}
 \Z[\chi_{1}^{\pm}, \chi_{2}^{\pm}] \rightarrow & \Z[\chi^{\pm}] \\
   \chi_{1} \mapsto & \chi \\
   \chi_{2} \mapsto & \chi^{n}.
\end{align*}

The following proposition is immediate from Proposition \ref{L:Edidin's tip}.

\begin{Proposition}\label{P:FnK}
The $\mc{T}$ equivariant $K$-theory of $\mb{F}_{n}$ is given by a $4$-tuple of elements in $(f_{1}, f_{2}, f_{3}, f_{4})$
where each $f_{i}$ (for $i = 1,2,3,4$) belongs to the ring $K_{*}(k)[ \chi^{\pm}]$, and satisfies the relations
\begin{multline*}
f_{1} -f_{2} = 0 \mod (1-\chi), \; f_{2} -f_{3} = 0 \mod (1- \chi^{2n}), \; f_{3} -f_{4} = 0 \mod (1- \chi),\; \\ \text { and } f_{4} - f_{1} = 0 \mod (1- \chi^{n}).
\end{multline*}
\end{Proposition}

The associated RS presentations are also easily computed. 
\begin{Proposition}[Reisner-Stanley presentation] \label{P:RSp}
The torus equivariant $K$-theory of wonderful rank one $\mc{G}$ compactifications have the following presentations.

\begin{itemize}
\item[] \[
 K_{\mc{T},*}(\mb{P}^{2}) =   \frac{K_{*}(k) [x_{\rho_{1}}^{\pm}, x_{\rho_{2}}^{\pm}, x_{\rho_{3}}^{\pm}]}{ ((x_{\rho_{1}}x_{\rho_{2}}x_{\rho_{3}}^{-2} -1), \; (x_{\rho_{1}}-1)(x_{\rho_{2}}-1)(x_{\rho_{3}}-1)) } \]
 
\item[] \[
K_{\mc{T},*}(\mb{P}^{1} \times \mb{P}^{1}) =  \frac{ K_{*}(k) [x_{\rho_{1}}^{\pm}, x_{\rho_{2}}^{\pm}, x_{\rho_{3}}^{\pm}, x_{\rho_{4}}^{\pm}]  } { ((x_{\rho_{1}}x_{\rho_{3}}^{-1} - x_{\rho_{2}}x_{\rho_{4}}^{-1}), \;(x_{\rho_{1}}-1)(x_{\rho_{3}}-1),\; (x_{\rho_{2}}-1)(x_{\rho_{4}}-1))} \]
 
\item[] \[
 K_{\mc{T},*}(\mb{F}_{n}) = \frac{K_{*}(k) [x_{\rho_{1}}^{\pm}, x_{\rho_{2}}^{\pm}, x_{\rho_{3}}^{\pm}, x_{\rho_{4}}^{\pm}]}{ ((x_{\rho_{1}}^{n}x_{\rho_{3}}^{-n} - x_{\rho_{2}}x_{\rho_{3}}^{n}x_{\rho_{4}}^{-1}), \;(x_{\rho_{1}}-1)(x_{\rho_{3}}-1), \;(x_{\rho_{2}}-1)(x_{\rho_{4}}-1))}.
\]

\end{itemize}  
\end{Proposition}

\begin{proof}
The proposition follows from the fact that we have a presentation $R(\mc{T})= R(T)/I$. 
As a result, by using Proposition \ref{L:Edidin's tip}, we get a presentation of the equivariant groups $K_{\mc{T},*}(\mc{X}) = K_{T,*}(\mc{X})/I$. 
The explicit description of the ideal $I$ then follows from the explicit calculations of the $R(T)$-module structure on the Reisner-Stanley presentation carried out in Sections \ref{sb:p2comp} - \ref{sb:fn}.

Let us work out the case of $\mb{P}^{2}$. The ideal is generated by $(1- \chi_{1}\chi_{2})$ and using the $R(T)$-module
structure of the RS presentation, worked out in  Section \ref{SS:p2},  we get \[R(\mc{T}) = \frac{\Z[\chi_{1}^{\pm}, \chi_{2}^{\pm}]}{(1- \chi_{1}\chi_{2})}.\]
The other cases are similar and this proves the proposition. 
\end{proof}   

\subsection{Weyl group action} 
The Weyl group of $\mc{G}$ acts on $\mc{X}$. Let $w_{0}$ denote the non trivial element in the Weyl group. It permutes
the torus fixed points and consequently acts on the ring $K_{*}(k)[\chi^{\pm}_{\sigma}]$ (recall the notation used in
Section \ref{SS:toric comp}). So, to understand the action of the Weyl group on the $\mc{T}$-equivariant $K$-theory we
need to understand the action on the torus fixed points.

\begin{Proposition}\label{P:weyl}
Let $w_{0}$ be the non-trivial element of the Weyl group of $\mc{G}$. Then the $w_{0}$ action on the compactifications
$\mc{X}$ are given as follows.
Let $\sigma_{ij}$ denote the maximal cones of the fans in Figures \ref{F:P2}, \ref{F:P1P1} and \ref{F:Fn} and $x_{\sigma_{ij}}$ denote the corresponding torus
fixed points of the toric varieties.

\begin{enumerate}
\item If $\mc{X} = \mb{P}^{1}$ then $w_{0}$ permutes the two torus fixed point. 
 
\item If $\mc{X}= \mb{P}^{2}$ then $w_{0}(x_{\sigma_{12}}) = x_{\sigma_{12}}$ and $w$ permutes the other two torus fixed points i.e. $ x_{\sigma_{23}} \stackrel{w_{0}}{\longleftrightarrow} x_{\sigma_{13}}.$

\item If $\mc{X}= \mb{P}^{1} \times \mb{P}^{1}$ then $x_{\sigma_{12}} \stackrel{w_{0}}{\longleftrightarrow} x_{\sigma_{34}}$ and it leaves the other two torus fixed points $x_{\sigma_{23}}$ and $x_{\sigma_{14}}$ invariant.

\item  If $\mc{X}= \mb{F}_{n}$ then $x_{\sigma_{12}} \stackrel{w_{0}}{\longleftrightarrow} x_{\sigma_{34}}$ and it leaves the other two torus fixed points $x_{\sigma_{23}}$ and $x_{\sigma_{14}}$ invariant.
    
\end{enumerate} 
\end{Proposition}   

\begin{proof}

The main ingredient in the proof is a ``dynamical'' interpretation of the torus fixed points in a toric variety. Let us
denote this torus fixed point, associated to a maximal cone $\sigma$ in the fan $\Delta$, by $x_{\sigma}$. Let $e \in X(\Delta)$ be any point in the dense torus $T$ of $X(\Delta)$. The point $x_{\sigma}$ is also the limit of the one parameter orbit $\nu(t) \cdot e$ where $\nu$ is any co-character in the interior of the cone $\sigma$.

The embedding $\iota:\mc{T} \rightarrow T$, considered in the diagram \eqref{E:embed}, defines a co-character
$i_{\chi} \defeq \iota(\chi)$ in $N_{T}$. We have a decomposition of lattices $N_{T} = \Z \cdot i_{\chi} \oplus \Z \cdot i_{\chi}^{\perp}$. We extend the action of $w_{0}$ to $N_{T}$ by defining $w_{0}(i) = -i$ and $w_{0}(i^{\perp})= i^{\perp}$. 
Any co-character of $\lambda \in N_{T^{b}}$ decomposes uniquely as $n\cdot i + m \cdot i^{\perp}$. Hence we can calculate the $w_{0}$ action by \[w_{0}(x_{\sigma}) = \lim_{t \rightarrow 0} w_{0}(\lambda)(t)\] where $\lambda$ is any co-character in the interior of $\sigma$.

The proposition follows by easy computations in a case-by-case analysis. 
\end{proof}

This allows us to completely calculate the $\mc{G}$-equivariant $K$-theory of rank one $\mc{G}$ wonderful varieties. 

\begin{Proposition}\label{P:PSL2} 
Let $\mc{X}$ be a rank one wonderful compactification of $\mc{G}$. Then its $\mc{G}$ equivariant $K$-theory is given as follows. 

\begin{enumerate}

\item If $\mc{X} = \mb{P}^{1}$ then $K_{\mc{G}, *}(\mc{X}) = K_{\ast}(k) \otimes R(T)$.

\item If $\mc{X} = \mb{P}^{2}$, the $K_{\mc{G}, *}(\mc{X})$ is given by  a collection of elements in $(f_{1}, f_{2})$ in the ring $\prod _{i=1}^{2} K_{*}(k)[ \chi^{\pm}]$ which satisfies the condition
$$
f_{1} -f_{2} = 0 \mod (1-\chi).
$$

\item If $\mc{X} = \mb{P}^{1} \times \mb{P}^{1}$, then $K_{\mc{G}, *}(\mc{X})$ is given by a collection of elements $(f_{1}, f_{2}, f_{3})$ in the ring $\prod _{i=1}^{3} K_{*}(k)[ \chi^{\pm}]$ which satisfies the condition
$$
f_{i} -f_{j} = 0 \mod (1 -\chi) \; \text { for all } i, j \in \{ 1,2,3\}.
$$

\item If $\mc{X} = \mb{F}_{n}$, then $K_{\mc{G}, *}(\mc{X})$ is given by a collection of elements in $(f_{1}, f_{2}, f_{3})$ in the ring $\prod _{i=1}^{3} K_{*}(k)[ \chi^{\pm}]$ which satisfies the conditions
$$
f_{1} -f_{2} = 0 \mod (1-\chi^{2n}), \text { and } f_{3} - f_{1} = 0 \mod (1-\chi^{n}).
$$
\end{enumerate}
\end{Proposition}

\begin{proof} The first assertion is well known and see Example \ref{E:gnthm} below for an outline of the proof.
The second and the third part of this proposition are direct consequences of Proposition \ref{P:weyl}, Proposition \ref{P:P1P1K} and Proposition \ref{P:P2K}. 

We consider the final part. 
Let $(f_{1}, f_{2}, f_{3}, f_{4})$ be an element in $\prod _{i=1}^{4} K_{*}(k)[ \chi^{\pm}]$ which satisfies the
conditions of Proposition \ref{P:weyl}. The invariants, described in Proposition \ref{P:weyl}, impose additional relations $ f_{1} = f_{3}$ and hence $f_{3} -f_{2} = f_{1} - f_{2} = 0 \mod (1-\chi^{2n}) $. This proves the proposition.    
\end{proof}

\begin{Remark}\label{R:rsw}
  It is tedious, but not difficult to work out the explicit Weyl group action on the Reisner-Stanley
  presentations. However the variables $x_{\rho}$ are not very well adapted to the this action. We will not use this action so we omit the details. 
\end{Remark}

\subsection{General case}
Let us now return to the general case of a reductive group $G$ and a fixed maximal torus $T \subset G$ and a $G$-variety
$X$. In this case the
$T$-equivariant $K$-theory admits a concrete description. The $G$-equivariant $K$-theory is much more subtle. The
problem is that $X$ admits a stratification by several $G$-orbits and only some of them contain $T$-fixed points.
In general the stabilizer of the Weyl-group action on the torus fixed points in each strata will be different. To the best of
our knowledge there is no uniform way to handle this issue; however see \cite{strickland2, bcmpaper} for the case of complete quadrics and
Section \ref{S:applications} for minimal rank symmetric varieties.
We have the following result.

\begin{Theorem}\label{T:main1}
Let $X$ be a smooth projective spherical $G$-variety. Let $T$ be a maximal torus of $G$. 
Then the set of $T$ fixed points $X^{T}$ is finite and the $T$ equivariant $K$-theory $K_{T,*}(X)$ is an ordered set of elements $(f_{x})_{x \in X^{T}}$ from $\prod_{x \in X^{T}} K_{*}(k) \otimes R(T)_{x}$ which are subject to the following additional congruences:
\begin{itemize}

\item $f_{x} - f_{y} = 0 \mod (1-\chi)$ when $x$ and $y$ are connected by a $T$-invariant curve of weight $\chi$. 

\item $f_{x} - f_{y} = f_{x} -f_{z}  = 0 \mod (1-\chi)$ and $f_{y} - f_{z} = 0 \mod (1-\chi^{2})$ where $\chi$ is a root
  of the pair $(G,T)$. The irreducible component of the subvariety $X^{\Ker(\chi)}$  which contains the points $x,y,z$
  is isomorphic to $\mb{P}^{2}$ and there is an element in the  Weyl group of $G$ that fixes $x$ and permutes the point $y$ and $z$.

\item $f_{x} - f_{y} = f_{y} -f_{z}  = f_{z} - f_{w} = f_{x} - f_{w} = 0 \mod (1-\chi)$ where $\chi$ is a root of the
  pair $(G,T)$. The irreducible component of the subvariety $X^{\Ker(\chi)}$ which contains the points $x,y,z,w$  is
  isomorphic to $\mb{P}^{1}\times \mb{P}^{1}$ and there is an element in the Weyl group of $G$ that fixes two points and permutes the other two. 

\item $f_{x} - f_{y} = f_{z} -f_{w} = 0 \mod (1-\chi)$ and $ f_{y} - f_{z} = 0 \mod (1-\chi^{2n})$ and $f_{x} - f_{w} =0
  \mod (1-\chi^{n})$, where $n \geq 1$, and $\chi$ is a root of the pair $(G,T)$. The irreducible component of the subvariety $X^{\Ker(\chi)}$  which contains the points $x,y,z,w$  is isomorphic to a ruled surface $\mb{F}_{n}$. There is an element in the Weyl group of $G$ that fixes the points $x$ and $w$ and permutes $z$ and $y$. 

\end{itemize}
The $G$-equivariant $K$-theory is given by the space of $W$-invariants in the $T$-equivariant $K$-groups. 
\end{Theorem}

\begin{proof}
Let $S$ denote any codimension one subtorus of $T$. 
We know that the irreducible components of $S$ fixed points $X^{S}$ are either  smooth curves $\mb{P}^{1}$ or one of the
compactifications listed in Proposition \ref{P:Brion-characterization}.
 
The formula, from Proposition \ref{P:desc}, 
\[ K_{T,*}(X^{S}) = K_{T/S, *}(X^{S}) \otimes_{R(T/S)} R(T)\]
reduces the problem to calculating $K_{T/S, \ast}(X^{S})$ which follows from Section \ref{S:Rank one case}. 
 
Note that in all these cases, the Weyl group of the neutral component of $C_{G}(S)$ (trivial when $C_{G}(S) = T$) embeds
as a subgroup of $W$. The proposition then follows from the Proposition \ref{P:weyl}.  
\end{proof}

\begin{Remark}
In the group compactification case, the connected components of $X^{S}$ are only smooth curves. In this case the
structure of the $T$-equivariant (and $G$-equivariant) $K$-theory has been worked out by Uma in \cite{Uma}.
\end{Remark}

\subsubsection*{Reisner-Stanley Presentation} \label{SS:RSpres}
Informally speaking, the Reisner-Stanley presentation is ``generated'' by patching the RS presentation of each
irreducible component of $X^{S}$ as $S$ varies over the codimension one subtori of $T$.  To formalize this, we construct 
a topological space $PL(X)$ associated to $X$. It is a two-dimensional topological space. We start a collection of points corresponding to
the $T$-fixed points of $X$. Then we glue the unit interval to pairs of points $\{x_{\alpha}, x_{\beta}\} \in X^{T}$ if
there is a $T$-stable curve passing through them. Next we consider subsets of three or four points in $X^{T}$. If
such a collection of points belong to an irreducible invariant surface then we glue a polyhedron, corresponding to the moment polytope
of the corresponding toric surface, with the vertices of the polyhedron glued to the torus fixed point. The boundary of
a two cell doesn't necessarily glue to one cells so the space $PL(X)$ is not a simplicial complex. The Weyl group
of $G$ acts on $PL(X)$. 

\begin{Proposition}[Reisner-Stanley presentation]\label{T:RSp}
Let $X$ be a smooth projective spherical $G$ variety and let $PL(X)$ be the topological space constructed in the
preceding paragraph. Let $s$ be an one or two dimensional cell of $PL(X)$ and we let $X_{s}$ denote the corresponding component. Then we have a map of rings \[ \varphi_{s}: K_{T/S, *}(X_{s}) \otimes_{R(T/S)} R(T) \rightarrow K_{T,*}(X^{T}).\]
The image of the colimit of the family of the maps $\varphi_{s}$ is identified with $K_{T,*}(X)$.
\end{Proposition}
 
\begin{proof}

  The $T/S$ invariant points of $X_{s}$ form a subset of $X^{T}$. Let us denote this subset by $X_{s}^{T/S}$.  The map
  $\varphi_{s}$ is the composition of maps, see Diagram (\ref{E:RS-map}), where the horizontal arrow is a result of Theorem \ref{VV}
  and the vertical arrow is inclusion. 
  \begin{equation}\label{E:RS-map}
    \begin{tikzcd}
       K_{T/S,*}(X_{s}) \otimes_{R(T/S)} R(T) \arrow[r] \arrow[rd, "\varphi_{s}"] & K_{T, \ast}(X_{s}^{T/S}) \otimes_{R(T/S)} R(T)
       \arrow[d]  \\
           & K_{T,*}(X^{T})
     \end{tikzcd}
   \end{equation}

The proposition now follows from the fact that an element $(f_{i}) \in K_{T,\ast}(X^{T})$ is in the image of
$K_{T,\ast}(X)$ when it satisfies the relations imposed by Theorem \ref{T:main1}. These relations only depend on the irreducible
component containing a given collection of points and the image of the map $\varphi_{s}$, from the corresponding
component, is a isomorphism (see Proposition \ref{P:RSp}). 
\end{proof}

\begin{Remark}
The previous proposition gives an immediate set of generators for $K_{T,*}(X)$. However in general the $W$ group
invariants are harder to extract because the variables used in this presentation are not well adapted to this group
action.
\end{Remark}

\begin{Example}\label{E:gnthm}
We consider the complete flag variety to illustrate the various presentations of the equivariant $K$-theory. 
We fix a connected reductive group $G$, a Borel subgroup $B$ and a maximal torus $T \subset B$. This fixes a root system $\Delta_{G}$ of $G$ and we denote the associated Weyl group by $W$. Given any simple root $\alpha \in \Delta_{G}$ we denote the associated element in $W$ by $s_{\alpha}$. We let $\mc{B}$ denote the complete flag variety $G/B$. 

We note that, using Proposition \ref{Faddeev-Shapiro} and the discussion after that we get $K_{G,*}(\mc{B}) = K_{B,*}(pt)$ and  $K_{B, *}(pt) = K_{T,*}(pt) = K_{*}(k) \otimes R(T)$.

\subsubsection*{GKM presentation} The $T$-fixed points and $T$-stable curves were described by Carrell in \cite [Theorem
F] {Carrell94}. It turns out, this data is given by the Bruhat graph $W_{G}$; the vertices are the $T$-fixed points of
$\mc{B}$ hence they are indexed by elements of the Weyl group. Two vertices $w$ and $w^{\prime}$ are joined by a
$T$-stable curve  if $w^{\prime} = s_{\alpha} \cdot w$ for some simple reflection $s_{\alpha}$ and $\ell(w)< \ell(w^{\prime})$. 

As a consequence of Theorem \ref{T:main1} we get 
\begin{equation*}
K_{T,*}(\mc{B}) = \left\{ (f_{w}) \in \prod_{w \in W} K_{*}(k) \otimes R(T) :
\begin{split} 
\mt{ where } f_{w} -f_{w^{\prime}} =0 \mod (1- \alpha) \mt{ whenever } \\
w^{\prime}= s_{\alpha} \cdot w \mt{ and } \ell(w) < \ell(w^{\prime}).
\end{split}
\right\}. 
\end{equation*}

The Weyl group  $W$ acts transitively on $\prod_{w \in W} K_{*}(k) \otimes R(T)$ by permuting the co-ordinates. So it suffices to check the congruences imposed at the identity. Any element in $K_{T,*}(\mc{B})^{W}$ is determined by choosing $f_{e} \in R(T)$. So, $K_{T,*}(\mc{B})^{W} = K_{*}(k) \otimes R(T)$. 

\subsubsection*{RS presentation} The topological space $PL(\mc{B})$ is just the topological realization of the
Bruhat graph. Let $s_{\alpha}$ be an edge of the Bruhat graph joining vertices $w_{0}$ and $w_{1} = w_{0} \cdot
s_{\alpha}$ (see Figure \ref{F:be} below). The component $X_{s_{\alpha}}$  is isomorphic to $\mb{P}^{1}$. We have RS
presentation, following Section \ref{SS:p1comp}, \[ K_{T/S,*}(X_{s_\alpha}) = \frac{K_{*}(k) [x_{\alpha_{1}}^{\pm},
    x_{\alpha_{2}}^{\pm}]}{ (x_{\alpha_{1}}-1)(x_{\alpha_{2}}-1) } = K_{\ast}(k) \]

The map $\varphi_{s_{\alpha}}$ in this case is just the identity on the corresponding end-points and the $T$-equivariant cohomology is identified
with \[ \prod_{e \in PL(\mc{B})} K_{\ast}(k) \otimes R(T)\]
 where the indexing set is the set of edges of the Bruhat graph. 

To compute the $G$-equivariant $K$-theory it suffices to look at the intersection of image of $\varphi$ and the sub-ring
$K_{*}(k) \otimes R(T)$ of $K_{T,*}(X^{T})$ corresponding to the identity $e \in W$. The Bruhat graph near $e$ is
depicted in Figure~\ref{F:brge}. Since all the edges emanate from $e$ we note that $G$-equivariant $K$-theory is
determined by the factor $K_{\ast}(k) \otimes R(T)$ at the identity. 
\begin{figure}[htp]
\begin{minipage}{0.45\textwidth}
\begin{center}
\begin{tikzpicture}[scale = 0.9, transform shape]
\begin{scope}

\path[draw, thick] (0,0) -- (2,0) ;
\draw (0,0) (0,0.4) node {$w_{0}$};
\filldraw (0,0) circle (2 pt);
\filldraw (2, 0) circle(2 pt);
\draw (2,0.4)  node {$w_{1}$};
\draw (1, -0.3) node {$ s_{\alpha}$};
\draw [dashed, gray, thick] (-1,0)--(0,0)  (-1,1)--(0,0)--(-1,-1);
\draw [dashed, gray, thick] (2,0)--(3,0)  (3,1)--(2,0)--(3,-1);
\end{scope}
\end{tikzpicture}
\caption{An edge of Bruhat graph connecting $w_{0}$ and $w_{1}= w_{0} \cdot s_{\alpha}$}
\label{F:be}
\end{center}
\end{minipage}
\hfill
\begin{minipage}{0.45\textwidth}
\begin{center}
\begin{tikzpicture}[scale = 0.8, transform shape]
\begin{scope}
\filldraw (0,0) circle(2pt); 
\draw (0.2,0.2)  node {$e$};
\draw [thick]  (0,0)--(0,2)  (0, 0)--(1.46, -1.46)  (0,0)--(-1.46, -1.46);
\draw (-1.8, -1.8) node{$\alpha_{1}$};
\filldraw[gray] (-1.46, -1.46) circle(2pt);
\draw (0, 2.4) node{$\alpha_{2}$};
\filldraw[gray] (0, 2) circle(2pt);
\draw (1.8, -1.8) node{$\alpha_{3}$};
\filldraw[gray] (1.46,-1.46) circle(2pt);
\draw [dashed, gray] (-1,3)--(0,2)  (1,3)--(0,2)  (2.46, -1.46)--(1.46,-1.46)  (1.46, -2.46)--(1.46,-1.46)   (-2.46, -1.46)--(-1.46,-1.46)  (-1.46, -2.46)--(-1.46,-1.46);
\end{scope}
\end{tikzpicture}
\caption{The local picture of Bruhat graph near the identity $e$. The vertices $\alpha_{i}$ are simple roots.}
\label{F:brge}
\end{center}
\end{minipage}
\end{figure}
\end{Example}

\begin{Remark} \label{S:Final}
  A consequence of Merkurjev's comparison result on ordinary and equivariant $K$-theory of algebraic varieties shows
  that smooth and projective varieties are \emph{equivariantly formal}\footnote{A notion first defined in
    \cite{GKM98}. }; see ~\cite[Proposition 3.6 and Theorem 4.2]{Merkurjev97}. In particular, we can recover the ordinary
  $K$-theory of $X$ from its $T$-equivariant $K$-theory by the formula
  \[ K_{\ast}(X) = K_{T,\ast}(X)/\mathfrak{m} \cdot K_{T, \ast}(X) \]
  where $\mathfrak{m}$ denotes the augmentation ideal of $R(T)$.
  
Alternatively, this also follows from Theorem 1.2 of \cite{JoshuaKrishna}.
\end{Remark}

%
%
%
%

\section{Application: Wonderful compactification of Symmetric Spaces of Minimal Rank}\label{S:applications}

In this section, we will apply our results to wonderful compactifications of symmetric spaces of minimal rank. 
We first recall some known facts about these varieties and refer the reader to the articles \cite{conspr, springerinvol,
  rich} for detailed proofs.

\subsection{Symmetric spaces and their compactifications}\label{S:symspth}
Throughout this section $G$ is a connected linear semisimple group of adjoint type and $\theta: G \rightarrow G$ be a
non-trivial involution. Let $H= G^{\theta}$ denote the subgroup of $\theta$-fixed points. The homogeneous space $G/H$ is
called a symmetric space. The group $H$ is reductive and without loss of generality we assume that it is connected
\footnote{Otherwise we replace $H$ by its connected component.}.

It turns out that the homogeneous space $G/H$ is a spherical $G$-variety and rank of $G/H$, as a spherical variety, is
always constrained by the inequality \[\rank(G) \geq \rank(H) + \rank(G/H).\]
A \emph{minimal rank} variety is a $G$-variety such that the above inequality is an equality. The rank doesn't change when one
passes to a spherical equivariant compactification of $G/H$. Throughout this section we will consider wonderful
compactifications of minimal rank symmetric spaces.  

\begin{Definition} A $\theta$-split torus of $G$ is a $\theta$-stable torus $S \subset G$ such that $\theta(s) =
  s^{-1}$. A $\theta$-split parabolic $P$ of $G$ is a parabolic subgroup such that $P \cap \theta(P)=$ Levi subgroup of $P$ (and $\theta(P)$).   
\end{Definition}

It turns out that for reductive groups nontrivial $\theta$-split torus always exist. Let us fix a maximal $\theta$-split
torus $S$ and a  minimal $\theta$-split parabolic $P$ containing $S$. Let $L$ and $P_{u}$ denote the Levi
component, and the unipotent radical of  $P$ respectively. It follows from construction that $L = C_{G}(S)$. The derived
group of $L$, denoted by $[L,L]$ is $\theta$-stable and it has no $\theta$-split torus. As a result, we conclude that $[L,L] \subset H$. We fix a maximal torus $T$ of $P$ containing $S$. As a consequence of the minimal rank assumption we see that $S \times (T \cap H) \rightarrow T$ is an isogeny.

Let $\Phi_{G}$ and $\Phi_{L}$ denote the root systems associated to the pairs $(G,T)$ and $(L,T)$. The positive roots,
the simple roots, and the associated Weyl groups are denoted by $\Phi^{+}_{G}$ (resp. $\Phi^{+}_{L}$), $\Delta_{G}$
(resp. $\Delta_{L}$) and $W_{G}$(resp. $W_{L}$) respectively. The involution $\theta$ acts on roots and we have a  subset \[\Phi^{-\theta} \defeq  \left\lbrace \alpha \in \Phi_{G}^{+} \;| \; \theta(\alpha) < 0 \right \rbrace .\] 
Since $[L,L] \subset H$ the action of $\theta$ is trivial on $\Phi_{L}$ and $W_{L}$. 
It turns out that we have a partition of positive roots \[ \Phi_{G}^{+} = \Phi^{- \theta} \cup \Phi_{L}^{+} \] and a
compatible partition of simple roots
\begin{equation}\label{E:root-partition}
  \Delta_{G} = \Delta_{L} \cup \Delta^{-\theta},
\end{equation}
where $\Delta^{-\theta} \defeq \Delta \cap \Phi^{-\theta}$.

The inclusion map $i: S \hookrightarrow T$ and the surjective homomorphism  $p: T \rightarrow S$ is defined by $p(t) = t\cdot \theta(t)^{-1}$ fits into a commutative diagram 
\begin{equation}\label{E:trs}
\begin{tikzcd}
 S \arrow{r}{i} \arrow[bend left]{rr}{*^{2}} & T \arrow{r}{p} & S 
\end{tikzcd}
\end{equation}
 where $*^{2}$ is the squaring map.
 
The injective map $p:M_{S} \rightarrow M_{T}$ on characters identifies $M_{S}$ with the subspace generated by the vectors $\alpha -\theta(\alpha)$. The nonzero vectors $i(\Phi_{G})$ in the image of the surjective map $i: M_{T} \rightarrow M_{S}$ generate a (possibly non-reduced) root system, denoted by $\Phi_{G/H}$, and the image $\Delta_{G/H} \defeq i(\Delta^{ - \theta})$ is a basis of this system. It turns out that the associated Weyl group $W_{G/H}$ of the root system equals $N_{H}(S)/C_{H}(S)$. The linear automorphisms of $\Phi_{G/H}$ are precisely the automorphisms of $M_{T}$ that preserve the subspace $p(M_{S})$. So we conclude that the 
$\theta$- invariant elements of the Weyl group $W_{G}$ surjects onto $W_{G/H}$. The partition of the root system, in
eqn.(\ref{E:root-partition}) above, produces a short exact sequence of groups 
\begin{equation}\label{E:sesq}
 1 \rightarrow W_{L} \rightarrow W_{G}^{\theta} \rightarrow W_{G/H} \rightarrow 1.
\end{equation}

Let us denote $T^{\theta} \defeq T \cap H$ which, under minimal rank assumption,  is a maximal torus of $H$. The results
of Brion and Joshua, see \cite[Section 1.4]{BrionJoshua08}, show that the possibly non-reduced root system $\Phi_{G/H}$ is
reduced, and the roots $\Phi_{H}$ is a subset of $\Phi_{G}$. The Weyl group $W_{H}$ equals the $\theta$-invariant
elements of $W_{G}$. So from the short exact
sequence eqn.(\ref{E:sesq}) leads to the short exact sequence
\begin{equation}\label{E:sesq-1}
 1 \rightarrow W_{L} \rightarrow W_{H} \rightarrow W_{G/H} \rightarrow 1.
\end{equation}

Let $X$ be the wonderful compactification of a minimal rank symmetric space $G/H$. The closure of the torus $S$ inside
$X$ is a toric variety, which we  will denote by $Y$. The fan associated to this toric variety is the
subdivision of $N_{S}$ (the space of $S$-co-characters) by the Weyl chambers of the root-system
$\Phi_{G/H}$. Let $Y_{0}$ denote the torus invariant affine open subset defined by the opposite Weyl chamber
\footnote{This Weyl chamber corresponds to the simple roots $-\Delta_{G/H}$.}. The Weyl group $W_{G/H}$ acts
transitively on the Weyl chambers so we have $Y = W_{G/K}\cdot Y_{0}$. The cone corresponding to the affine subset
$Y_{0}$ is of maximal dimension and we denote the unique torus fixed point associated to this cone, in $Y$, by
$z_{0}$. 

There is a complete classification of irreducible symmetric spaces and their compactifications. In each case the
topological space $PL(X)$, introduced in Section \ref{SS:RSpres}, is a graph. More precisely, all the codimension-one torus
stable components are only smooth curves. A complete parametrization of the torus fixed points of $X$ and the torus
stable curves connecting these points in known. 

\begin{Proposition}[See Lemma 2.1.1 \cite{BrionJoshua08}]\label{L:BrionJoshua}
We continue with the notation above. The torus fixed points and the curves of wonderful symmetric space $X$ and the toric variety $Y$ admit following parametrization:
\begin{enumerate}[(i)]

\item The $T$-fixed points of $X$ are exactly the points $w \cdot z_0$ where 
$w\in W_G/ W_{L}$. The torus fixed points of $Y$ are exactly $w \cdot z_0$ where  $w \in W_H/W_L = W_{G/H}$.  

\item For any positive root $\alpha \in \Phi^{+}_{G} \setminus \Phi^{+}_{L}$, there exists unique irreducible $T$-stable 
curve $C_{\alpha \cdot z_0}$ connecting $z_0$ and $\alpha_{z_{0}}$. The torus $T$ acts on $C_{\alpha \cdot z_0}$
via the character $\alpha$. The curve is isomorphic to $\mb{P}^{1}$ and we call these curves \emph{Type 1} curves.

\item For any simple root $\gamma =\alpha - \theta(\alpha) \in \Delta_{G/H}$, there exists unique irreducible 
$T$-stable curve $C_{\gamma \cdot z_{0}}$ connecting $z_0$ and $s_\alpha s_{\theta(\alpha)}\cdot z_0$. 
The torus $T$ acts on $C_{\gamma \cdot z_{0}}$ by the character $\gamma$. 
The curve is isomorphic to $\mb{P}^{1}$ and we call these curves \emph{Type 2} curves.

\item The irreducible $T$-stable curves in $X$ are precisely the $W_G$-translates of the curves
$C_{\alpha \cdot z_0}$ and $C_{\gamma \cdot z_0}$. They are all isomorphic to $\PP^1$.

\item The irreducible $T$-stable curves in $Y$ are the $W_{G/H}$-translates of the curves
$C_{\gamma \cdot z_0}$.

\end{enumerate}

\end{Proposition}

\subsection{Equivariant $K$-theory of symmetric spaces}\label{SS:Kthsym}
In this section we will explore the structure of equivariant $K$-theory of wonderful compactifications of
symmetric-spaces. Our main tools will be the previous proposition and Theorem \ref{T:main1}.

\subsubsection{T-equivariant $K$-theory}\label{SS:TEQK}

We start with an immediate consequence of Proposition \ref{L:BrionJoshua}.    

\begin{Corollary}\label{C:step1}
The $T$-equivariant $K$-theory $K_{T,*}(X)$ of $X$ is isomorphic to the space of 
tuples $(f_{w\cdot z_0}) \in \prod_{w \in W_G/W_L} K_{*}(k) \otimes R(T)$ such that 
\[
 f_{w \cdot z_0} - f_{w^{\prime}\cdot z_0} = 
  \begin{cases} 
   0 \mod (1- \alpha) & \text{if } w^{-1}w^{\prime}= s_{\alpha} \\
   0 \mod (1- \alpha \cdot \theta(\alpha)^{-1})    & \text{if }   w^{-1}w^{\prime}= (s_\alpha \cdot s_{\theta(\alpha)})^{\pm}
  \end{cases} .
\]

The $T$-equivariant $K$-theory $K_{T,*}(Y)$ of the toric variety $Y$ is isomorphic to the space of 
tuples $(f_{w\cdot z_0}) \in \prod_{w \in W_H/W_L} K_{*}(k) \otimes R(T)$ such that 

\[
 f_{w \cdot z_0} - f_{w^{\prime}\cdot z_0} = 
     0  \mod (1- \alpha \cdot \theta(\alpha)^{-1})     \text{ and }   w^{-1}w^{\prime}= (s_\alpha \cdot s_{\theta(\alpha)})^{\pm}.
 \]
\end{Corollary}

Next, we relate the $T$-equivariant $K$-theories of $X$ and $Y$. Let $W^{H}$ denote the minimal coset representatives of
$W_{G}/W_{H}$ (recall that $W_{H} = W_{G}^{\theta}$). We have the following proposition. 

\begin{Proposition}\label{P:Tkthsym}
There is an isomorphism of rings 
\[ \prod_{w \in W^{H}} K_{T,*}(Y) \cong K_{T,*}(X) \]
which is compatible with the $K_{*}(k) \otimes R(T)$-module structure on both sides.   
\end{Proposition} 
 
\begin{proof}
Consider the chain of Coxeter groups $ W_{L} \subset W_{H} \subset W_{G}$. It introduces a map on the quotient spaces
$\pi: W_{G}/W_{L} \rightarrow W_{G}/W_{H}$. We identify $W_{G}/W_{L}$ (resp. $W_{G}/W_{H}$) set-theoretically with
$W^{L}$ (resp. $W^{H}$) - the minimal length coset representatives. This identification defines a section, denoted by
$\gamma$, to the map $\pi: W^{L} \rightarrow W^{H}$. Note that $W^{L}$ is a $W_{H}/W_{L} = W_{G/H}$ torsor over $W^{H}$.
    
Let $\varphi$ denote the composition of ring homomorphisms in eqn.(\ref{E:spttr}) where the first map is canonical and
the second map is the rearrangement map 
\begin{equation}\label{E:spttr}
\begin{tikzcd}
K_{T,*}(X) \arrow[rr, "\varphi", bend left] \arrow[r] & K_{T,*}(X^{T}) \arrow{r} & \prod_{u \in W^{H}}(\prod_{v \in \pi^{-1}(u)} K_{*}(k) \otimes R(T)).
\end{tikzcd}
\end{equation}
We denote the image of $\varphi$ inside  as the ring $R_{\varphi}$, and for any $u \in W^{H}$ we denote the intersection
$R_{\varphi} \cap (\prod_{v \in \pi^{-1}(u)}K_{\ast} \otimes R(T))$ by $R_{\varphi}^{u}$. The $W$-action on the torus
fixed points $X^{T}$ translates into a $W$-action on the ring $\prod_{u \in W^{H}}(\prod_{v \in \pi^{-1}(u)} K_{*}(k)
\otimes R(T))$ and, using Proposition \ref{L:BrionJoshua} and Corollary \ref{C:step1}, the action introduces ring isomorphisms
\begin{equation} \label{E:ringiso}
  \gamma_{w} : R_{\varphi}^{u} \rightarrow R^{u}_{\varphi} \rightarrow R^{u \cdot w}_{\varphi}.
\end{equation}
We use the right action here because the subgroup $W_{L} \subset
W$ must act trivially.

The description of $T$-equivariant $K$-groups, in Corollary \ref{C:step1}, show that $K_{T,\ast}(Y) = R_{\varphi}^{e}$ where $e \in W^{H}$ is the element of smallest length (i.e., it corresponds to the coset of $[H]$).
So we get an inclusion map
\begin{equation} \label{E:incltoric}
 i: K_{T,\ast}(Y) \rightarrow R_{\varphi}. 
\end{equation}

Given any $w \in W^{H}$ we define the map $i_{w}$ by the composition $\gamma_{w} \circ i$ and this defines the 
map
\[ i_{W}: \prod_{w \in W^{H}}K_{T, \ast}(Y) \rightarrow R_{\varphi}\]
in each co-ordinate.

The map $i_{W}$ is clearly injective and it is surjective by $W$-equivariance. This proves the proposition. 
\end{proof}

\subsubsection{Alternate description in terms of simplicial complex}
The toric variety $Y$ is a compactification of the maximal anisotropic torus $S$ and as a $T$-variety the torus $T/S$ acts trivially on
$Y$. Applying Proposition \ref{P:desc} in this case, we get  
\begin{equation}\label{E:teqt}
K_{T,*}(Y) = K_{S,*}(Y) \otimes R(T/S).
\end{equation}

In \cite{BDP90}, the authors associate a simplicial-complex $\ms{C}_{Y}$ to a smooth toric variety $Y$ (see Definition 5
of \cite{BDP90}). The
complex $\ms{C}_{Y}$ encodes the geometry of the fan (which defines the toric variety $Y$). The setting in \cite{BDP90}
is that of equivariant cohomology but their argument is geometric and it works for $K$-theory as well.

One can associate a Reisner-Stanley algebra to a simplicial complex, purely combinatorially over any coefficient ring, and the algebra admits a
direct-sum decomposition into submodules where the summands are combinatorially defined. Using the coefficient ring $K_{\ast}(k)$ we get
decomposition \[K_{S,*}(Y) = \oplus_{\Delta \in
    \ms{C}_{Y}} K_{S,*}(Y)_{\Delta}\] where the summand $K_{S,*}(Y)_{\Delta}$ is a $K_{\ast}(k)$-module consists of the
monomials which are supported the simplex $\Delta \subset \mc{C}_{Y}$.

This decomposition can now be exploited for the toric variety $Y$ that appears in the description of the
wonderful-compactification of the symmetric space $X$. 

\begin{Proposition}
  The $T$-equivariant $K$-theory of $X$ admits the following direct sum- decomposition
  \begin{equation}\label{E:tdirs}
K_{T,*}(X) = \prod_{W_{H}} \; \left (\bigoplus_{\Delta \in \ms{C}_{Y}}  K_{S,*}(Y)_{\Delta} \otimes R(T/S) \right ). 
\end{equation}
where $\ms{C}_{Y}$ is the simplex associated to the toric variety $Y$. 
\end{Proposition}

\begin{Remark}
In the group case this recovers Lemma 2.8 of \cite{Uma}.
\end{Remark}


\subsubsection{$G$-equivariant $K$-theory}
We consider the $G$-equivariant $K$-theory and we will provide two descriptions of it. The first one is a direct
consequence of Proposition \ref{P:Tkthsym}. 
 
\begin{Proposition}\label{C:W_H invariants}
The $G$-equivariant $K$-theory $K_{G,*}(X)$ of $X$ 
is isomorphic to $W_H$-invariants of the $T$-equivariant $K$-theory of the toric variety $Y$.
\end{Proposition}

\begin{proof}
The $G$-equivariant $K$-theory of $X$ is given by the formula  $K_{G,*}(X) = K_{T,*}(X)^{W_{G}}$. We have a
set-theoretic splitting of $W_{G} = W^{H} \times W_{H}$.

It is clear that any element of $K_{T,\ast}(Y)^{W_{H}}$ embeds into $K_{G,\ast}(X)$. Conversely, it follows from
Proposition \ref{P:Tkthsym}, that after translating by an element of $W^{H}$ any element of $K_{G, \ast}(X)$ must embed
into $K_{T, \ast}(Y)$. Further it must also be invariant with respect to $W_{H}$-actions. So the injection
$K_{T,\ast}(Y)^{W_{H}}$ is also a surjection.

\end{proof}

\begin{Remark}
  One can also prove the above Proposition directly from the description of the torus stable curves and points on $X$
  and $Y$ outlined in Proposition \ref{L:BrionJoshua}.
\end{Remark}

\subsubsection{A refined  description of $G$-equivariant $K$-theory}
It turns out that one can refine Proposition \ref{C:W_H invariants} even further. The Weyl Group $W_{G/H}$ acts
transitively on the Weyl chambers of the root system $\Phi_{G/H}$ and as a result the toric variety $Y$ admits a cover
by the $W_{G/H}$- translates of an affine open set $Y_{0}$ (recall $Y_{0}$ was is the affine-set corresponding to the
anti-dominant Weyl chamber).

The subgroup $W_{L} \subset W_{H}$ acts trivially on $Y$ (because it acts trivially on the fan) and hence the action of
$W_{H}$ on $Y$ factors via $W_{H}/W_{L} = W_{G/H}$. The following equality then follows easily \[
 K_{S,*}(Y)^{W_{H}} = K_{S,*}(Y)^{W_{G/H}} = K_{S,*}(Y_{0}).\]
This leads to the following structure theorem for $G$-equivariant $K$-theory. 

\begin{Proposition}\label{P:structure}
The $G$-equivariant $K$-theory of $X$ is 
\[ K_{G,*}(X) = K_{S,*}(Y_{0}) \otimes R(T/S)^{W_{H}}. \]
\end{Proposition}

\begin{proof}
In the light of Proposition \ref{C:W_H invariants}, it suffices to show that \[K_{T,*}(Y)^{W_{H}} = K_{S,*}(Y_{0}) \otimes R(T/S)^{W_{H}}.\]

The key-idea here is that, as a consequence of minimal-rank condition, we can identify $T/S$ with the maximal torus of
$H$. As a result, there is a Steinberg basis of $R(T/S)$ over $R(T/S)^{W_{H}}$. We will denote the elements of the
Steinberg basis by $\left \lbrace e_{w} \right \rbrace_{w \in W_{H}}$. 
  
It follows from eqn.(\ref{E:teqt}) that  $K_{S,*}(Y_{0}) \otimes R(T/S)^{W_{H}} \subset K_{T,*}(Y)^{W_{H}}$ and we
will show the other inclusion. The Steinberg basis forms a basis $\left \lbrace 1 \otimes e_{w} \right \rbrace_{w \in W_{H}}$ forms a basis of
$K_{T,*}(Y)$ over $K_{S,*}(Y) \otimes R(T/S)^{W_{H}}$. As noted above $W_{H}$ acts transitively on  $K_{S,*}(Y)$ so
taking $W_{H}$ invariants we get $K_{S,*}(Y_{0}) \otimes R(T/S)^{W_{H}}$.
\end{proof}

\subsubsection{Multiplicative structure constants}
As noted in Proposition \ref{P:Tkthsym}, the factor $K_{T,*}(Y)$ essentially determines the $T$-equivariant $K$-theory of $K_{T,\ast}(X)$. The toric variety $Y$
is determined by the reduced and irreducible root system $\Phi_{G/H}$ with Weyl group $W_{G/H}$.
The root system $\Phi_{G/H}$ determines an adjoint algebraic group $\Gamma(\Phi_{G/H})$. Let $B \Gamma(\Phi_{G/H})$ denote the complete flag variety of $\Gamma(\Phi_{G/H})$. Then, using \cite{kly},  we can identify the toric variety $Y$ with the closure of the general torus orbit of the maximal torus in the flag variety $B\Gamma(\Phi_{G/H})$. The following lemma is immediate.

\begin{Lemma}
Let $S$ denote the maximal torus of $\Gamma(\Phi_{G/H})$ compatible with the root system $\Phi_{G/H}$. Then $K_{S,*}(B\Gamma(\Phi_{G/H})) = K_{S,*}(Y)$. 
\end{Lemma}

The torus equivariant $K$-theory of flag variety has deep combinatorial structures. In particular Kostant and Kumar, in
\cite{kok}, show that the $K$-theory $K_{S,0}(B\Gamma(\Phi_{G/H}))$ admits a remarkable basis, called the Schubert
basis, over $R(S)$. Let us denote this basis by $\{ [ \mc{O}_{w}] \}_{w \in W_{G/H}}$; it is indexed by the Weyl group $W_{G/H}$. 

\begin{Proposition}
  The torus equivariant $K$-theory, $K_{T,\ast}(X)$, admits a natural Schubert-basis over $R(T)$. 
\end{Proposition}

\begin{proof}
  We note that eqn.(\ref{E:teqt}) expresses the fact that $K_{T, \ast}(X)$ is obtained by extension of scalars from
  $K_{S,\ast}(Y)$. The existence of Schubert basis for $K_{S,\ast}(Y)$ shows that we have a presentation
  \[
    K_{S,\ast}(Y) = R(S)[\{\mc{O}_{w}\}_{w \in W_{G/H}}] .
  \]
  The proposition then follows immediately.
\end{proof}

\begin{Remark}
The Schubert basis (in the equivariant setting)  exhibits \emph{positivity} phenomenon  (see
\cite{gr} for precise definitions and conjectures) and it has deep combinatorial structure. The previous proposition
shows that in the case of wonderful compactifications of minimal rank symmetric varieties, one can recover the
structure constants from that of a lower dimensional flag-variety.  
\end{Remark}

In the group case, when $X$ is the wonderful compactification of $G\times G/\Delta(G)$, the Schubert basis
$ \{ \mc{O}_{w}\}$ corresponds to the Schubert-basis of the flag-variety of $G$.

%
%
%
%

\appendix
\section{ Equivariant Intersection Theory} \label{SS:cecg}

In this section we will show that, for a $G$-variety $X$, the torus equivariant $K$-theory
determines the torus equivariant Chow-theory. In contrast to the rest of the paper, in this section by $K$-theory we
mean the equivariant Grothendieck group with rational
coefficients, i.e,  $K_{T,\Q}(X) \defeq K_{T,0} \otimes \Q$. By Chow-theory we mean the equivariant Chow-ring with
rational coefficients; as defined by Graham and Edidin. Our
arguments are rather formal in nature and it works in a more general setting\footnote{For example for certain non-smooth
$T$-skeletal  varieties one may substitute operational theories.}. The precise hypothesis on the space $X$, which is always
satisfied by smooth spherical projective varieties, are explained below. Our main tool is the equivariant Riemann-Roch
map defined by Edidin and Graham, see \cite{edgr2}.

\subsubsection*{Notation}
We consider a $G$-variety $X$, and fix a maximal torus $T$ of $G$ with
character lattice $M$. For example, we have $ K_{T,\Q}(pt) \defeq \Q[M]$ whereas
$ A^{*}_{T,\Q}(pt) \defeq \Sym_{\Z}(M) \otimes \Q$.

Let $I \subset K_{T,\Q}(pt)$ (resp. $J \subset A^{*}_{G, \Q}(X)$) denote the corresponding augmentation ideals. So $I$
(resp. $J$) is generated by $1- \chi$ (resp. $\chi$) for characters $\chi \in M$.  Let $\widehat{K_{T,\Q}(X)}$
denote the completion of $K_{T, \Q}(X)$ with respect to the $I$-adic filtration. We view
$\prod_{i=0}^{\infty}A^{i}_{T,\Q}(X)$ as a completion of the Chow-theory $\oplus_{i=0}^{\infty}A^{i}_{T,\Q}(X)$ with
respect to the sub-modules $A^{[n]}_{T,\Q}(X) \defeq \prod_{i=n}^{\infty}A^{i}_{T, \Q}(X)$. The $T$-equivariant Riemann-Roch map,
denoted by $\tau^{T}$,  maps  $\tau^{T}: K_{T, \Q}(X) \rightarrow \prod_{i=0}^{\infty}A^{i}_{T,\Q}(X)$ and it induces an isomorphism, also denoted by $\tau^{T}$, between the completions
\begin{equation} \label{E:RIEMANN-ROCH}
  \tau^{T}:\widehat{K_{T,\Q}(X)} \rightarrow \prod_{i=0}^{\infty}A^{i}_{T,\Q}(X). \end{equation}

In particular, when $X=pt$,  we identify $\prod_{i=0}^{\infty}A^{i}_{T,\Q}(pt)$ with the ring of formal power-series $\sum_{i=0}^{\infty} a_{i} \chi^{i}$ (where $\chi \in M$) and the map $\tau^{T}$ evaluated at the element $1- \chi$ is given by 
\begin{equation} \label{E:RRmap}
\tau^{T}(1-\chi) = \sum_{i=1}^{\infty} (-1)^{i+1} \chi^{i}/ (i+1)! . 
\end{equation}

\subsubsection*{Assumption}
The precise assumptions on the nature of the space $X$ are as follows.

\begin{itemize}

\item We assume that the locus of $T$ fixed points of $X$ is finite and the natural restriction maps $K_{T, \Q}(X)
  \rightarrow K_{T,\Q}(X^{T})$ and $A^{*}_{T,\Q} \rightarrow A^{*}_{T,\Q}(X^{T})$ are injective. When $X$ is a smooth
  spherical $G$-variety these conditions are always satisfied.

\item The embedding $K_{T, \Q}(X) \rightarrow K_{T,\Q}(X^{T})$ is defined by finitely many congruence conditions of the
  form $f_{i} = f_{j} \mod p_{ij}$ where $f_{i}, f_{j}, p_{ij}$ are elements of $K_{T,\Q}(pt)$, and $p_{ij}$ belongs to the augmentation ideal $I$. 
\end{itemize}

When the equivariant $K$-theory (resp. Chow-theory) of $X$ satisfies the second condition above we say that the $K$-theory
(or Chow-theory) is \emph{commensurable}. We will assume that the above assumptions are always satisfied.

\begin{Lemma}\label{L:INJ}
  The following assertions are true for equivariant $K$-theory as well as equivaraint Chow-theory.
  \begin{itemize}
\item The $I$-adic filtration on $K_{T, \Q}(X)$ and the filtration induced on $K_{T,\Q}(X)$ as a sub-module of
  $K_{T,\Q}(X^{T})$ with its $I$-adic filtration are equivalent.

\item The completion map $K_{T,\Q}(X) \rightarrow \widehat{K_{T,\Q}(X)}$ is injective.
\end{itemize}


\end{Lemma}

\begin{proof}
The $K_{T,\Q}(pt)$ module $K_{T,\Q}(X^{T})$ is finitely generated and $K_{T,\Q}(pt)$ is a noetherian ring. So the first
assertion is a consequence of the Artin-Rees lemma and the second assertion is a consequence of the Krull Intersection
theorem (see \cite[Chapter 8]{matsumura}).
The proof for equivariant Chow-theory is verbatim.  
\end{proof}

The following proposition is the main result of this section. 
\begin{Proposition}
Suppose the equivariant $K$-theory of $X$ is commensurable then so is the equivariant Chow-theory of $X$. A commensurable
presentation of the equivariant $K$-theory in terms of finitely many congruence conditions defines a
commensurable presentation of the equivariant Chow-theory. 
\end{Proposition}

\begin{proof}
Let us assume that the equivariant $K$-theory of $X$ is commensurable and it is defined by finitely many congruence conditions of the form $f_{i} = f_{j} \mod p_{ij}$ (where $p_{ij}$ belongs to the augmentation ideal $I$) with a finite index set $(i,j) \in \mc{S}$. 

In the Diagram \eqref{E:KA} below, the first and the third squares commute because completion is functorial. The second
square commutes because the equivariant Riemann-Roch map is functorial. The maps $i$ (resp. $j$) are injective
by assumption, and the maps $\widehat{i}$ (resp. $\widehat{j}$) are injective by Lemma \ref{L:INJ}.

\begin{equation}
\label{E:KA}
\begin{tikzcd}[row sep=normal, column sep=large]
K_{T,\Q}(X) \arrow{d}{i} \arrow{r}  & \widehat{K_{T,\Q} (X)} \arrow{r}{\tau^{T}} \arrow{d}{\widehat{i}} &  \prod_{i=0}^{\infty} A^{i}_{T,\Q}(X)  \arrow{d}{\widehat{j}} &  A^{*}_{T,\Q}(X) \arrow{l} \arrow{d}{j} \\
K_{T,\Q}(X^{T}) \arrow{r} & \widehat{K_{T,\Q}(X^{T})} \arrow{r}{\tau^{T}} & \prod_{i=0}^{\infty}A^{i}_{T,\Q}(X^{T})  &  A^{*}_{T,\Q}(X^{T}) \arrow{l}  
\end{tikzcd}
\end{equation}

For any non-negative integer $n$, let $\mc{F}^{n}(K_{T,\Q}(X))$ denote the sub-module $I^{n} \otimes K_{T,\Q}(X)$ and
$\mc{F}^{n}(\widehat{K_{T,\Q}(X)})$ the corresponding sub-module in $ \widehat{K_{T,\Q}(X)}$. We identify the quotients $\widehat{K_{T, \Q}(X)}/\mc{F}^{n}(\widehat{K_{T,\Q}(X)})$ with  $K_{T,\Q}(X)/\mc{F}^{n}(K_{T,\Q}(X))$ and similarly $\prod _{i=0}^{\infty}A^{i}_{T, \Q}(X)/ A^{[n+1]}_{T, \Q}(X) =
\oplus_{i}^{n} A^{i}_{T,\Q}(X)$. By construction, the equivariant Riemann-Roch isomorphism (see
eqn.(\ref{E:RIEMANN-ROCH})) is continuous with respect to the completions so we get a commutative square

\begin{equation} \label{E:RRdim}
\begin{tikzcd}
K_{T,\Q}(X)/\mc{F}^{n}(K_{T,\Q}(X)) \arrow[r, "i_{n}"] \arrow[d, "\tau_{n,m}^{T}"] & K_{T,\Q}(X^{T})/\mc{F}^{n}(K_{T,\Q}(X^{T})) \arrow[d,"\tau^{T}_{n,m}"] \\
\oplus_{i=0}^{m} A^{i}_{T,\Q}(X) \arrow{r}{j_{m}} & \oplus_{i=0}^{m} A^{i}_{T,\Q}(X^{T}) 
\end{tikzcd}
\end{equation}
where $m-n$ is bounded for $m \geq n \gg 0$.

The terms in the right-hand column  of the Diagram (\ref{E:RRdim}) are free modules, and the
image of $i_{n}$ stabilizes as $n \gg 0$ because the $K$-theory presentation is commensurable. As a result
the image of $j_{m}$ also stabilizes as $m \gg 0$.

 The associated-graded ring of $\prod_{i=0}^{\infty}
A^{i}_{T,\Q}(X^{T})$ (with respect to the filtration $A^{[n]}_{T,\Q}(X^{T})$) is given by $A^{\ast}_{T,\Q}(X^{T})$. As a
result we note that the image of $A^{\ast}_{T,\Q}(X)$ is
determined inside $A^{\ast}_{T,\Q}(X)$ by finitely many relations. Moreover, passing to the associated
graded ring the image of the finitely many relations \[ \{ \tau_{n,m}^{T}(f_{i}) = \tau_{n,m}^{T}(f_{j}) \mod \tau_{n,m}^{T}(p_{ij})\}_{(i,j) \in
\mc{S}}\] determine a commensurable presentation of the Chow-ring. 
\end{proof}  

\begin{Remark}
We can recover Brion's calculation of equivariant Chow-theory of smooth projective spherical variety using Theorem
\ref{T:main1} and the above proposition. 

A similar result also holds for $G$-equivariant theories because the $G$-equivariant theory is determined by the
invariants of the geometric action of the Weyl-group on the torus fixed points.  
\end{Remark}

\bibliographystyle{alpha}
\bibliography{EKT.bib}

\end{document}